%% file: main.tex
\documentclass[onefignum,onetabnum,reqno]{siamart190516}
\usepackage[normalem]{ulem}
\usepackage{geometry}
\usepackage{amsmath}
\usepackage{physics}
\renewcommand{\sech}{\mathrm{sech}}
\input{Misc/shared_siam}
\begin{document}
\maketitle
\begin{abstract}
We propose a multilevel tensor-train (TT) framework for solving nonlinear partial differential equations (PDEs) in a global space–time formulation. While space–time TT solvers have demonstrated significant potential for compressed high-dimensional simulations, the literature contains few systematic comparisons with classical time-stepping methods, limited error convergence analyses, and little quantitative assessment of the impact of TT rounding on numerical accuracy. Likewise, existing studies fail to demonstrate performance across a diverse set of PDEs and parameter ranges. In practice, monolithic Newton iterations may stagnate or fail to converge in strongly nonlinear, stiff, or advection-dominated regimes, where poor initial guesses and severely ill-conditioned space–time Jacobians hinder robust convergence. We overcome this limitation by introducing a coarse-to-fine multilevel strategy fully embedded within the TT format. Each level refines both spatial and temporal resolutions while transferring the TT solution through low-rank prolongation operators, providing robust initializations for successive Newton solves. Residuals, Jacobians, and transfer operators are represented directly in TT and solved with the adaptive-rank DMRG algorithm. The proposed method preserves the compression and scalability of TT representations while substantially improving convergence and stability.  Numerical experiments for a selection of nonlinear PDEs — including Fisher–KPP, viscous Burgers, sine–Gordon, and KdV — cover diffusive, convective, and dispersive dynamics, demonstrating that the multilevel TT approach consistently converges where single-level space–time Newton iterations fail. 
In dynamic, advection-dominated (nonlinear) scenarios, multilevel TT surpasses single-level TT, achieving high accuracy with significantly reduced computational cost, specifically when high-fidelity numerical simulation is required.
\end{abstract}

\begin{keywords}
  Quantized Tensor train, space-time solver, nonlinear partial differential equations
\end{keywords}

 \begin{AMS}
   68Q12, 68Q09 
   ,65F55 
   ,65M06 
   ,14N07 
  ,65M30 
\end{AMS}

\section{Introduction}
The numerical solution of nonlinear time-dependent PDEs remains a significant challenge, particularly for problems exhibiting sharp gradients, wave phenomena, or stiff dynamics. In this work, we present and rigorously validate a novel computational framework that significantly advances the state-of-the-art for such problems. We compare three distinct numerical approaches: the {classical time-stepping (CT)} method, which serves as our baseline; a {single-level tensor-train space-time (SL-TT)} method, representing a direct low-rank advancement of monolithic schemes; and our primary contribution, the {multi-level tensor-train space-time (ML-TT)} method.

The mathematical foundation of our approach relies on tensor train (TT) decompositions, which provide a powerful framework for mitigating the curse of dimensionality in high-dimensional problems. The TT format was introduced as a computationally efficient alternative to canonical and Tucker decompositions, expressing a $d$-dimensional tensor as a chain of low-rank cores \cite{oseledets2011tensor}. Subsequent developments included the TT-cross approximation for black-box tensor approximation \cite{Oseledets-Tyrtyshnikov:2010} and the Quantized TT (QTT) format, which achieves logarithmic compression by folding long vectors into higher-order tensors \cite{kazeev2012low}. These advances have enabled efficient solutions of high-dimensional PDEs \cite{Dolgov:2012,khoromskij2011tensor} and have been integrated with iterative solvers like DMRG \cite{Schollwock:2011} and AMEN \cite{dolgov2014alternating} for linear systems in TT format.

The monolithic space–time approach formulates the evolution problem as a single coupled system over the entire space–time domain, offering an alternative to sequential time stepping. By solving for all time levels simultaneously, the method avoids the classical accumulation of local truncation errors associated with step-wise integration and supports a uniformly implicit treatment of nonlinearities over the entire time horizon. When combined with tensor-train (TT) representations, the space–time formulation can exploit separability in both spatial and temporal dimensions, yielding substantial compression of the space–time degrees of freedom and reducing the effective dimensionality of the problem. Previous works have demonstrated their potential, including space-time TT approaches for Fokker--Planck equations \cite{Dolgov:2012}, nonlinear methods for differential Riccati equations \cite{Breiten:2021}, and space-time methods for advection-diffusion and reaction-diffusion systems via spectral/collocation and tensor train methods \cite{Adak:2024, Adak:2025, Arenstein:2025}. Separately, implicit time marching on TT manifolds with stepwise truncation has been developed to control rank growth during integration \cite{Rodgers:2023}.

However, existing single-level TT space-time methods often face key limitations. Their performance and robustness can be highly sensitive to the initial guess for stiff nonlinear problems. Furthermore, there are few direct comparisons with classical time-stepping across diverse nonlinear PDEs; systematic error studies that separate discretization error from TT rounding effects are limited; and monolithic Newton iterations may lack robustness in stiff or advection-dominated regimes. A notable gap is the performance of these methods in challenging {nonlinear hyperbolic regimes}. For instance, while the viscous Burgers equation has been studied in its parabolic regime, a thorough quantitative evaluation of monolithic TT solvers in its hyperbolic regime—where nonlinear advection dominates, leading to shock formation—has not been previously presented. In such regimes, the resolution of hyperbolic problems may benefit from block-time-domain decomposition and the use of the Parallel in Time Stable Bi-Conjugate Gradient algorithm \cite{riahi2022pitsbicg} to mitigate the ill-posedness arising from high bond dimensions and further accelerate the solution through parallel computing. Furthermore, recent quantum-inspired developments have introduced a space-time solver that achieves logarithmic runtime scaling and enables accurate data-driven predictions via algorithms like MPS-DMD \cite{peddinti2025quantum}. These approaches highlight the effectiveness of global space-time representations in capturing complex spatio-temporal correlations with minimal degrees of freedom.

Our work directly addresses these gaps. While the single-level TT method (SL-TT) is a contribution in its own right, our primary advancement is the {multi-level TT (ML-TT) framework}. This method leverages a hierarchy of coarse-to-fine space-time grids to systematically generate high-quality initial guesses, dramatically improving the robustness and convergence speed of the Newton solver. 

Our study goes beyond proof of concept, demonstrating validated solutions across a broad spectrum of nonlinear dynamics: parabolic behavior (c.f., Fisher-KPP~\cite{Fisher:1937} and viscous Burgers in the diffusion-dominated regime), hyperbolic behavior (viscous Burgers in the advection-dominated, shock-forming regime and the sine-Gordon equation~\cite{Bratsos:2008}), and dispersive wave propagation (c.f., the Korteweg-de Vries equation~\cite{Zabusky:1965}). 
For each case, we perform detailed mesh refinement studies and analyze the impact of low-rank truncation. The ML-TT method consistently outperforms both the classical and single-level approaches, showing particular strength in stabilizing convergence for demanding problems. Its successful application to a broad class of nonlinear problems constitutes a key result of this paper and demonstrates the method's capability to handle strong nonlinearities and sharp gradients where existing monolithic low-rank solvers have not been thoroughly tested.

At first glance, the proposed multilevel strategy may appear related to classical space–time multigrid methods~\cite{Horton:1995,RapakaS:2020}. However, the present approach is fundamentally different in both purpose and construction. Classical multigrid methods operate on the linearized operator and rely on residual restriction, coarse-grid correction, and cycling strategies to accelerate convergence of a linear solver. In contrast, our multilevel framework is applied exclusively at the nonlinear level and is used solely to generate high-quality initial guesses for monolithic Newton iterations on successively refined space–time grids. No residual transfer, restriction to coarse-grid, or multilevel cycling is performed, and each level is solved independently to convergence, after which the solution is prolongated to the successive finer grid to serve as an initial guess for the Newton method. This distinction is crucial in the tensor-train setting, where the dominant failure mode is nonlinear stagnation or divergence rather than slow linear convergence.

The monolithic formulation allows for a flexible choice of time integrators. Our TT-formulation includes both the first-order implicit Euler scheme for its robustness and the second-order temporal schemes, such as Crank-Nicolson and Newmark methods, where higher accuracy is required to resolve wave dynamics. The global systems arising from all these schemes are seamlessly incorporated into the tensor-train framework.


The paper is organized as follows: Section~\ref{sec:tt} provides a necessary background on tensor train formats. Section~\ref{sec:TT_space_time} details the monolithic space-time discretization, the Newton solver, and the novel multi-level algorithm. Section~\ref{sec:results} presents our extensive numerical experiments, comparing all three methods across the different PDE classes and regimes. Finally, Section~\ref{sec:conclusion} offers concluding remarks.

\section{Tensor train (TT) methods}\label{sec:tt}
The tensor train (TT) format is a separable representation for high-dimensional arrays that mitigates the curse of dimensionality by expressing a $d$-way tensor as a chain of small 3D “cores.” Let $\mathcal{A}\in\mathbb{R}^{n_1\times\cdots\times n_d}$ with entries $\mathcal{A}(i_1,\ldots,i_d)$. A TT representation writes
\begin{equation}
\mathcal{A}(i_1,\ldots,i_d)
=\sum_{\alpha_1=1}^{r_1}\cdots\sum_{\alpha_{d-1}=1}^{r_{d-1}}
G^{(1)}(1,i_1,\alpha_1)\,
G^{(2)}(\alpha_1,i_2,\alpha_2)\cdots
G^{(d)}(\alpha_{d-1},i_d,1),
\label{eq:tt}
\end{equation}
where the TT-\emph{cores} are $G^{(k)}\in\mathbb{R}^{r_{k-1}\times n_k\times r_k}$ with $r_0=r_d=1$. The integers $r_k$ are the \emph{TT-ranks} (also called \emph{bond dimensions}); small ranks imply strong separability and yield substantial compression. 
Contracting adjacent cores along a rank index—the \emph{rank--core product}, denoted $\bowtie$—is the basic building block of TT arithmetic. For two consecutive cores,
\[
\big(G^{(k)} \bowtie G^{(k+1)}\big)(\alpha_{k-1},i_k,i_{k+1},\alpha_{k+1})
=\sum_{\alpha_k=1}^{r_k} G^{(k)}(\alpha_{k-1},i_k,\alpha_k)\,G^{(k+1)}(\alpha_k,i_{k+1},\alpha_{k+1}).
\]
Thus $\bowtie$ contracts the \emph{bond} of size $r_k$. The bond dimension $\chi=\max_k r_k$ controls both expressivity (admissible correlations) and computational cost.
With the \emph{rank--core product} $\bowtie$,
\[
\mathcal{A}(i_1,\ldots,i_d)
=G^{(1)}(:,i_1,:)\ \bowtie\ G^{(2)}(:,i_2,:)\ \bowtie\ \cdots\ \bowtie\ G^{(d)}(:,i_d,:).
\]
Matrices $A\in\mathbb{R}^{N\times M}$ with $N=\prod n_k$ and $M=\prod m_k$ use \emph{TT-matrix} (a.k.a.\ matrix product operator) cores $M^{(k)}\in\mathbb{R}^{r_{k-1}\times n_k\times m_k\times r_k}$:
\[
A\big((i_1,\ldots,i_d),(j_1,\ldots,j_d)\big)
=\sum_{\alpha_1,\ldots,\alpha_{d-1}}
\prod_{k=1}^d M^{(k)}(\alpha_{k-1},i_k,j_k,\alpha_k).
\]
A TT-matrix times a TT-vector remains in TT with ranks (before rounding) multiplying componentwise: $r_k^{Ax}\le r_k^{A}\,r_k^{x}$.

\subsection{Basic algebra in TT}
Algebraic operations in the tensor-train (TT) format can be expressed through
core-wise combinations of TT chains.
Let
\[
\mathcal{X} = G^{(1)} \bowtie \cdots \bowtie G^{(d)}, 
\qquad
\mathcal{Y} = H^{(1)} \bowtie \cdots \bowtie H^{(p)},
\qquad
A = M^{(1)} \bowtie \cdots \bowtie M^{(d)}
\]
be tensors in TT format with TT ranks $r_k^{\mathcal{X}}$, $r_k^{\mathcal{Y}}$, and $r_k^{A}$,
respectively, where
$r_0^{\mathcal{X}} = r_d^{\mathcal{X}} = r_0^{\mathcal{Y}} = r_p^{\mathcal{Y}} = r_0^{A} = r_d^{A} = 1$.
Here, $\mathcal{X}$ and $\mathcal{Y}$ are TT vectors with a single physical index per core,
whereas $A$ is a TT matrix (also referred to as a matrix product operator) whose cores
$M^{(k)}$ carry two physical indices corresponding to the input and output modes.
For addition, elementwise multiplication, and matrix--vector or matrix--matrix products,
the tensors must have compatible physical dimensions, i.e., identical mode sizes and $d=p$.
For the Kronecker product, the physical dimensions may differ and $d$ need not equal $p$.
\begin{table}[h]
\centering
\begin{tabular}{@{}p{4.5cm}p{5.5cm}p{4.2cm}@{}}
\toprule
Operation & Core construction & Rank behavior \\ \midrule
Addition ($\mathcal{X} + \mathcal{Y}$)
& Block-diagonal concatenation of cores
& $r_k = r_k^{\mathcal{X}} + r_k^{\mathcal{Y}}$; TT-rounding required \\

Elementwise (Hadamard) product ($\mathcal{X} \odot \mathcal{Y}$)
& Rank-wise Kronecker coupling of cores with shared physical indices
& $r_k \le r_k^{\mathcal{X}} r_k^{\mathcal{Y}}$; TT-rounding required \\

Matrix--vector ($A\mathcal{X}$)/ matrix--matrix product
& Local contraction over physical indices
& $r_k \le r_k^{A} r_k^{\mathcal{X}}$; TT-rounding required \\

Kronecker product ($\mathcal{X} \otimes \mathcal{Y}$)
& Concatenation of TT cores via a rank-one bond
& No rank growth; no rounding required \\ \bottomrule
\end{tabular}
\caption{Basic algebraic operations in the tensor-train (TT) format.}
\label{tab:tt_algebra}
\end{table}

Table~\ref{tab:tt_algebra} summarizes the main TT algebraic operations and their rank
behavior. Addition results in additive rank growth, while elementwise and matrix-based
products introduce multiplicative rank growth through Kronecker-type or contraction
operations. In all such cases, TT-rounding is required to control the ranks.
In contrast, the Kronecker product admits an exact TT representation obtained by
concatenating the TT cores of the operands through a trivial rank-one bond.
The operations listed in Table~\ref{tab:tt_algebra} are defined as follows.
\begin{align*}
\mathcal{X} + \mathcal{Y}
&=
\big[ G^{(1)} \;\; H^{(1)} \big]
\;\bowtie\;
\Big[
\begin{matrix}
G^{(2)} & 0 \\
0 & H^{(2)}
\end{matrix}
\Big]
\;\bowtie\;\cdots\;\bowtie\;
\Big[
\begin{matrix}
G^{(d)} \\
H^{(d)}
\end{matrix}
\Big],
\\
\mathcal{X} \odot \mathcal{Y}
&=
\big( G^{(1)} \otimes_{\mathrm{r}} H^{(1)} \big)
\;\bowtie\;\cdots\;\bowtie\;
\big( G^{(d)} \otimes_{\mathrm{r}} H^{(d)} \big),
\\
A\mathcal{X}
&=
\big( M^{(1)} \star G^{(1)} \big)
\;\bowtie\;\cdots\;\bowtie\;
\big( M^{(d)} \star G^{(d)} \big),
\\
\mathcal{X} \otimes \mathcal{Y}
&=
G^{(1)} \;\bowtie\;\cdots\;\bowtie\; G^{(d)}
\;\bowtie\;
H^{(1)} \;\bowtie\;\cdots\;\bowtie\; H^{(p)}.
\end{align*}
Here, $\otimes_{\mathrm{r}}$ denotes a Kronecker product taken over the TT rank indices and
applied independently for each physical index, while the core-wise contraction product
$\star$ represents a local contraction over the shared physical index between TT cores,
with the associated TT rank indices combined multiplicatively; both operators are defined
explicitly below, and the same construction applies to TT-matrix $\times$ TT-matrix
multiplication.
\begin{align*}
\big(G^{(k)} \otimes_{\mathrm{r}} H^{(k)}\big)
(\alpha_{k-1}\beta_{k-1},\, i_k,\, \alpha_k\beta_k)
&=
G^{(k)}(\alpha_{k-1}, i_k, \alpha_k)\,
H^{(k)}(\beta_{k-1}, i_k, \beta_k).
\\
\big(M^{(k)} \star G^{(k)}\big)
(\alpha_{k-1}\beta_{k-1},\, i_k,\, \alpha_k\beta_k)
&=
\sum_{j_k=1}^{n_k}
M^{(k)}(\alpha_{k-1}, i_k, j_k, \alpha_k)\,
G^{(k)}(\beta_{k-1}, j_k, \beta_k).
\end{align*}

\subsection{Rounding and memory compression}
The storage of a TT tensor is
\[
\text{nnz}(\text{cores})=\sum_{k=1}^d r_{k-1}\,n_k\,r_k
\;\approx\; \mathcal{O}(d\,n\,r^2)
\quad\text{(if } n_k\equiv n,\ r_k\equiv r).
\]
This replaces the exponential $\mathcal{O}(n^d)$ cost by nearly linear scaling in $d$ when ranks remain moderate. After algebraic operations, \emph{TT-rounding} (a sequence of QR+SVD truncations) reduces ranks to meet a tolerance $\varepsilon$, with work $\mathcal{O}(d\,n\,r^3)$.

\subsection{Quantized tensor train (QTT)}
QTT further compresses long vectors and matrices by \emph{folding} large physical indices into many small (often binary) modes before applying TT. For a vector of length $N$, choose a base $b$ (typically $b=2$) and write $N=b^{q}$; index $i\in\{0,\ldots,N-1\}$ has base-$b$ digits $i=(\iota_1,\ldots,\iota_q)$ with $\iota_k\in\{0,\ldots,b-1\}$. The folding map $\mathcal{F}_b:\mathbb{R}^{N}\to\mathbb{R}^{b\times\cdots\times b}$ ($q$ modes) reshapes $x$ to $\tilde x(\iota_1,\ldots,\iota_q)=x(i)$. Applying TT to $\tilde x$ yields a QTT with cores of size $r_{k-1}\times b\times r_k$ and storage
\[
\mathcal{O}\big(q\,b\,r^2\big)=\mathcal{O}\big((\log_b N)\,b\,r^2\big),
\]
which is \emph{logarithmic} in the original length $N$ when ranks remain bounded. For matrices $A\in\mathbb{R}^{N\times M}$ with $N=b^{q}$ and $M=b^{p}$, fold both row and column indices to obtain a QTT-matrix with $q+p$ small modes (each of size $b$). Many structured vectors and operators (e.g., exponentials, sinusoids, Toeplitz/Laplacian-like operators) admit small QTT ranks, enabling dramatic compression and fast arithmetic. All basic operations (addition, Hadamard, matvec) proceed as in TT, with temporary rank growth controlled by QTT-rounding.
 
The $\mathcal{O}(d\,n\,r^2)$ storage for TT—and $\mathcal{O}((\log_b N)\,r^2)$ for QTT with small base $b$—enables robust low-rank exploitation and memory compression, making TT/QTT practical tools for high-dimensional PDEs and related problems.
\section{Multi-level TT-space-time method}\label{sec:TT_space_time}
We consider a general nonlinear initial boundary value problem for a scalar field \( u = u(\mathbf{x}, t) \), where \( \mathbf{x} \in \Omega \subset \mathbb{R}^d \), \( t \in (0,T] \):

\begin{align}
\frac{\partial^m u}{\partial t^m} + \mathcal{L} u + \mathcal{N}(u) &= s(\mathbf{x}, t), && \mathbf{x} \in \Omega,\ t \in (0,T], \label{eq:pde} \\
\mathcal{B}(u, \nabla u, \ldots) &= g(\mathbf{x}, t), && \mathbf{x} \in \partial\Omega,\ t \in (0,T], \label{eq:bc} \\
u^j(\mathbf{x}, 0) &= u^j_0(\mathbf{x}), && \mathbf{x} \in \Omega, \quad j=0,1,\hdots,m-1,\label{eq:ic}
\end{align}
where, $\mathbf{x} = (x_1, x_2, \ldots, x_d)^\top \in \mathbb{R}^d $ is the spatial coordinate, $\mathcal{L}$ is a linear spatial differential operator,
$\mathcal{N}(u)$ is a nonlinear spatial differential operator acting on $u$, $s$ is a source term, $\mathcal{B}$ is a boundary operator (Dirichlet, Neumann, Robin, or periodic), and $u_0(\mathbf{x})$ is the initial condition. In this work we consider cases with $m=1,2$.

To simplify the description, we consider a one-dimensional space $\mathbf{x}=x\in [x_a, x_b]$, and it is straightforward to extend the present methodology to multiple space dimensions. Let us discretize the space and time domains into $N_x\equiv 2^{q_x}$ and $N_t=2^{q_t}$ equally spaced cells or intervals with a spacing of $\Delta x=(x_b-x_a)/N_x$ and $\Delta t= T/N_t$, respectively, and let $q=q_x+q_t$. Let $x_i=x_a+i\Delta x/2, i=1,2,\hdots,N_x$ be the center of each cell and $t_n=n\Delta t, n=1,2,\hdots,N_t$ be the end of each time interval. Let us define $u_i^n \coloneqq u(x_i,t_n)$ and denote the spatial solution vector at time $t_n$ as $U^n = [u_1^n,u_2^n,\hdots,u_{N_x}^n]^\top\in \mathbb{R}^{N_x}$, and the space-time solution vector as $\Uvec = [U^{1\top}, U^{2\top}, \hdots, U^{N_t\top}]^\top\in \mathbb{R}^{N_t\cdot  N_x}$. Similarly, let $s_i^n \coloneqq s(x_i,t_n)$ and denote the spatial source vector at time $t_n$ as $S^n = [s_1^n,s_2^n,\hdots,s_{N_x}^n]^\top\in \mathbb{R}^{N_x}$, and the space-time source vector as $\Svec = [S^{1\top}, S^{2\top}, \hdots, S^{N_t\top}]^\top\in \mathbb{R}^{N_t\cdot N_x}$.
Let $\Ld$ and $\Nd$ be the discretized linear and nonlinear operators, and define 
\begin{align*}
    F(U,S) \coloneqq \Ld U+\Nd(U)-S.
\end{align*}

To facilitate the boundary conditions, we refer to ghost cells which lie outside the domain (with indices 0 on the left and $N_x+1$ on the right) and have the same width as the interior cells.
The ghost cell values $ u_0$ and $ u_{N+1}$ can be expressed in terms of the values at the interior cells $ u_1$ and $ u_N$, along with the boundary condition (denoted with subscript $b$). 
\begin{align}
 u_0 &= a_1 u_1 + b_1 ,\qquad 
 u_{N+1} = a_2 u_{N} + b_2 ,\\
a_1, a_2 &= 
\begin{cases}
	-1 & \text{Dirichlet BC} \\
    +1 & \text{Neumann BC} \\
    0 & \text{no BC}
\end{cases}, \quad 
b_1, b_2 = 
\begin{cases}
	2 u_b & \text{Dirichlet BC} \\
    \left.\delux{u}{n}\right|_b \Delta n & \text{Neumann BC} \\
    0 & \text{no BC}
\end{cases}, \label{def:aLbL}
\end{align}
where $a_1$ and $a_2$ depend on the type of boundary conditions at left and right boundaries; $b_1$ and $b_2$ depend on the inhomogeneous boundary values at left and right boundaries.
Here, $\delux{u}{n} = \nabla u \cdot \hat{n}$,  $\Delta n$ and $\hat{n}$ are the mesh spacing and the unit vector along the outward normal direction, respectively. 

\subsection{All-at-once space-time discretization settings}
The space-time discretization uses the following vectors and matrices: $\mathbf{e}_1 = [1,\ 0,\ \dots,\ 0]^\top \in \mathbb{R}^{N_t}, \mathbf{e}_2 = [0,\ 1,\ 0,\ \dots,\ 0]^\top \in \mathbb{R}^{N_t}, I_x \in \mathbb{R}^{N_x\times N_x}, I_t \in \mathbb{R}^{N_t\times N_t}$ are the identity matrices,  $D_t\in \mathbb{R}^{N_t\times N_t}$ is a lower bidiagonal Toeplitz matrix with diagonal entries \(1\) and subdiagonal entries \(-1\), $D_{tt}\in \mathbb{R}^{N_t\times N_t}$ is a lower tridiagonal Toeplitz matrix with diagonal entries \(1\), first subdiagonal entries \(-2\) and second subdiagonal entries \(1\), $J_t\in \mathbb{R}^{N_t\times N_t}$ is a lower bidiagonal Toeplitz matrix with diagonal entries \( \frac{1}{2}\) and subdiagonal entries \( \frac{1}{2}\), and $K_t\in \mathbb{R}^{N_t\times N_t}$ is a lower tridiagonal Toeplitz matrix with diagonal entries \( \frac{1}{4}\), first subdiagonal entries \( \frac{1}{2}\) and second subdiagonal entries \( \frac{1}{4}\). Matrix and TT representations of all the operators are presented in Appendix \ref{sec:tt_operators}.

All discretizations that follow yield the nonlinear system for the space-time solution vector $\Uvec$, solved with Newton's method:
\begin{align}
        f(\Uvec) = A(\Uvec) + B \Uvec - \mathcal{C} = 0 \label{eqn:f_U}
\end{align}
where $A(\Uvec)$ is nonlinear, $B$ is linear, and $\mathcal{C}$ aggregates the initial conditions, source term, and any contributions from inhomogeneous boundary condtions; these components are defined separately for each case.

For $m=2$, the discrete approximation of $\delsux{u}{t}$ requires $U^{(-1)}$, the spatial vector at $t=-\Delta t$, which is numerically approximated about $t=0$ via a truncated Taylor expansion to keep $U_{tt}$ consistent with the temporal discretization's accuracy. For a second-order disretization of $\delsux{u}{t}$, $U^{(-1)}$ is approximated as:
\begin{equation}
    U^{(-1)} = U^0 - \Delta t \,U_t^0 + \frac{\Delta t^2}{2} \,U_{tt}^0 - \frac{\Delta t^3}{6} U_{ttt}^0 + \mathcal{O}(\Delta t^4), \label{eqn:U-1}
\end{equation}
where, the time derivatives $U_{tt}^0$ and $U_{ttt}^0$ are evaluated by differentiating Eq. \eqref{eq:pde} and substituting the initial conditions $U^0$ and $U^0_t$. 

\subsubsection{First order in time discretization}
For $m=1$, implicit Euler discretization of Eq. \eqref{eq:pde} yields:
\begin{align*}
    \frac{U^{n}-U^{n-1}}{\Delta t} + F(U^{n}, S^{n}) = 0, \quad n=1,2,\hdots,N_t.
\end{align*}
The full space-time solution vector $\Uvec$ is obtained from the system:
\begin{align}
    \frac{1}{ \Delta t} \left[(D_t\otimes I_x) \Uvec - \mathbf{e}_1 \otimes U^0\right] + (I_t \otimes \Ld) \Uvec + (I_t \otimes \Nd) \left( \Uvec \right) &= \Svec.\nonumber
\end{align}
Rearrange the above equation in the form of Eq. \eqref{eqn:f_U}: 
\begin{align}
    \underbrace{\Delta t (I_t \otimes \Nd)}_{A} \left( \Uvec \right) + \underbrace{(D_t\otimes I_x + \Delta t I_t \otimes \Ld)}_{B} \Uvec - \underbrace{(\mathbf{e}_1 \otimes U^0 + \Delta t \Svec)}_{\mathcal{C}} &= 0. \label{eqn:Euler_m1}
\end{align}

For $m=2$, Eq. \eqref{eqn:f_U} includes two initial conditions, $U^0$ and $U^0_t$, and includes the second derivative operator $D_{tt}$:
\begin{align}
    \underbrace{\Delta t^2 (I_t \otimes \Nd)}_{A} \left( \Uvec \right) + \underbrace{(D_{tt}\otimes I_x + \Delta t^2 I_t \otimes \Ld)}_{B} \Uvec - \underbrace{(\mathbf{e}_1 \otimes (2U^0 - U^{(-1)}) - \mathbf{e}_2 \otimes U^0 + \Delta t^2 \Svec)}_{\mathcal{C}} &= 0, \label{eqn:Euler_m2}
\end{align}
where, $U^{(-1)} = U^0 - \Delta t U_t^0 + \frac{\Delta t^2}{2} U_{tt}^0$, with $U_{tt}^0$ evaluated from Eq. \eqref{eq:pde} using the initial data $U^0$.
\subsubsection{Second order in time discretization}
For $m=1$, Eq. \eqref{eq:pde} is discretized with Crank-Nicolson scheme:
\begin{align*}
    \frac{U^{n}-U^{n-1}}{\Delta t} + \frac{1}{2}\left[F(U^{n},S^{n}) + F(U^{n-1},S^{n-1})\right] = 0,  \quad n=1,2,\hdots,N_t.
\end{align*}
The full space-time solution vector $\Uvec$ satisfies the nonlinear system:
\begin{align*}
    \frac{1}{ \Delta t} \left[(D_t\otimes I_x) \Uvec - \mathbf{e}_1 \otimes U^0\right] + (J_t \otimes \Ld) \Uvec + (J_t \otimes \Nd) \left( \Uvec \right) - (J_t\otimes I_x) \Svec + \frac{1}{2}(\mathbf{e}_1 \otimes F(U^0,S^0)) = 0.
\end{align*}
Rearrange the above equation in the form of Eq. \eqref{eqn:f_U}: 
\begin{align}
    \underbrace{\Delta t (J_t \otimes \Nd)}_{A} \left( \Uvec \right) + \underbrace{(D_t\otimes I_x + \Delta t J_t \otimes \Ld)}_{B} \Uvec - \underbrace{\left(\mathbf{e}_1 \otimes \left[U^0-\frac{\Delta t}{2}F(U^0,S^0)\right] + \Delta t (J_t\otimes I_x) \Svec\right)}_{\mathcal{C}} = 0. \label{eqn:CN_m1}
\end{align}

For $m=2$, Eq. \eqref{eq:pde} is discretized with an implicit, second-order, unconditionally stable Newmark scheme \cite{Newmark:1959}: 
\begin{align*}
    \frac{U^{n}-2U^{n-1}+U^{n-2}}{\Delta t^2} &+ \frac{1}{4}\left[F(U^{n},S^{n}) + 2F(U^{n-1},S^{n-1}) + F(U^{n-2},S^{n-2})\right] = 0,  ~~ n=1,2,\hdots,N_t. \label{eqn:NM_discretization}
\end{align*}
For $n=1$, $U^{(-1)}$, the spatial solution vector at $t=-\Delta t$, is approximated using Eq. \eqref{eqn:U-1}.
The full space-time solution vector $\Uvec$ then satisfies:
\begin{align}
    &\underbrace{\Delta t^2 (K_t \otimes \Nd)}_{A} \left( \Uvec \right) + \underbrace{(D_{tt}\otimes I_x + \Delta t^2 K_t \otimes \Ld)}_{B} \Uvec - \underbrace{\left(\mathbf{e}_1 \otimes C^1 - \mathbf{e}_2 \otimes C^2+ \Delta t^2 (K_t\otimes I_x) \Svec\right)}_{\mathcal{C}} = 0, \label{eqn:NM_m2}\\
    &C^1 = \left[2U^0 - U^{(-1)} -\frac{\Delta t^2}{2} F(U^0,S^0) -\frac{\Delta t^2}{4} F(U^{(-1)},S^{(-1)})\right], \quad C^2 = \left[U^0 + \frac{\Delta t^2}{4} F(U^0,S^0)\right].\nonumber
\end{align}

\subsection{Newton method}
We use Newton's method to solve Eq. \eqref{eqn:f_U} iteratively. In each iteration $k$, the solution is updated as follows:
\begin{align}
    \Uvec_{k+1} &= \Uvec_k + \omega_k \Uvec'_k,\label{eqn:Uvec_kp1}\\
    J_k \Uvec'_k &= -  f(\Uvec_k), \label{eqn:dU_Newton}
\end{align}
where, $\omega_k \in (0,1]$ is a relaxation factor and the Jacobian matrix $J$ is given by,
\begin{align}
    J \equiv \delux{f(\Uvec)}{\Uvec} = A'(\Uvec) + B. 
    \label{eqn:Jacobian}
\end{align}
Alternatively, we may solve:
\begin{align}
    J_k \Uvec_{k+1} &= J_k \Uvec_k - \omega_k f(\Uvec_k) := \mathcal{R}_k.\label{eqn:U_Newton}
\end{align}
For line search method, the above equation is reformulated to avoid repeated solves for varying $\omega_k$:
\begin{align}
    J_k \widetilde{\Uvec}_{k+1} &= J_k \Uvec_k - f(\Uvec_k) := \mathcal{R}_k = A'(\Uvec_k) \Uvec_k - A(\Uvec_k) + \mathcal{C} \label{eqn:U_Newton_line_solve}
    \\
    \Uvec_{k+1} &= (1-\omega_k) \Uvec_k + \omega_k \widetilde{\Uvec}_{k+1}
    \label{eqn:U_Newton_line_update}
\end{align}
Although Eqs.~\eqref{eqn:Uvec_kp1}–\eqref{eqn:dU_Newton} and
Eqs.~\eqref{eqn:U_Newton_line_solve}–\eqref{eqn:U_Newton_line_update}
are mathematically equivalent, we prefer the latter formulation.
In this form, the quantity $\mathcal{R}_k$ appears as a modified right-hand side in the Newton correction equation, as it eliminates the linear term $B\Uvec_k$ and retains only the nonlinear contribution together with the initial and boundary terms. In our numerical experiments, the absence of the linear $B\Uvec_k$ term is observed to significantly improve the conditioning of the Newton correction and the nonlinear convergence behavior. In addition, this reformulation consistently results in lower TT ranks for $\mathcal{R}_k$ compared to the full right-hand side
$f(\Uvec_k)$, leading to a corresponding reduction in computational cost.

For the solution of linear system in tensor train format (cf. Eq. \eqref{eqn:dU_Newton} or Eq. \eqref{eqn:U_Newton_line_solve}), we employ a two-site density matrix renormalization group (DMRG) method \cite{Oseledets:2012} (also known as the Modified Alternating Linear Scheme, MALS \cite{holtz2012alternating}) to solve linear systems in tensor train format. The method is rank-adaptive and automatically determines the TT-ranks. If needed, the solution ranks can be truncated according to input parameters such as the rounding error $\eps_{\mathrm{TT}}$ and the maximum TT-rank $\chi$. Additionally, within each two-site DMRG sweep, the global Newton system
\(
J_k \Uvec_{k+1} = b
\)
is reduced to a sequence of effective local linear systems of the form
\begin{equation}
A_{\mathrm{loc}} \, x = b_{\mathrm{loc}},
\end{equation}
where \( A_{\mathrm{loc}} \) and $b_{\mathrm{loc}}$ are the projected Jacobian $J_k$ and right hand side vector $b$ associated with the active TT cores.
For strongly nonlinear or advection-dominated problems, \( A_{\mathrm{loc}} \) may become ill-conditioned or nearly rank-deficient, leading to stagnation or instability of the Newton iteration.

To enhance robustness, we employ Tikhonov regularization at the level of these local systems.
Specifically, instead of solving the unregularized least-squares problem, we compute
\begin{equation}
x = \arg\min_{y} \left\{ \| A_{\mathrm{loc}} y - b_{\mathrm{loc}} \|_2^2
+ \alpha \, \| y \|_2^2 \right\},
\end{equation}
where \( \alpha > 0 \) is a regularization parameter.
Using the singular value decomposition \( A_{\mathrm{loc}} = U \Sigma V^{\top} \), the regularized solution can be written explicitly as
\begin{equation}
x = V \left( \Sigma^2 + \alpha I \right)^{-1} \Sigma U^{\top} b_{\mathrm{loc}}
= A_{\mathrm{loc},\alpha}^{\dagger} \, b_{\mathrm{loc}},
\end{equation}
which corresponds to a Tikhonov-regularized Moore--Penrose pseudoinverse.
This regularization damps the contribution of small singular values of \( A_{\mathrm{loc}} \), stabilizing the local solves while preserving the low-rank structure of the TT representation. 

Algorithms \ref{algo:Newton_left} and  \ref{algo:Newton_right} outline the steps to solve Eq. \eqref{eqn:f_U} in light of these settings.

\noindent

\begin{algorithm}[H]
\caption{Newton method for solving Eq.~\eqref{eqn:f_U} using Eq. \eqref{eqn:dU_Newton}}
\label{algo:Newton_left}
\begin{algorithmic}[1]
    \State \textbf{Input:} Initial guess $\Uvec_0$, TT-rounding inputs $\eps_{\mathrm{TT}}, \chi$ Newton solver inputs $n_{iter}, \eps_{res}, \eps_{cor}$, line search parameters $s, n_{line}$, linear system solver inputs $\eps_{\mathrm{DMRG}}, n_{sweeps}$, Tikhonov regularization parameter $\alpha$
    \State \textbf{Output:} Solution $\Uvec_{k+1}$

    \State Set $\eps_{\mathrm{TT},k}=10^{-3}$ 
    \State \textbf{Step 1:} Compute initial Residue $f(\Uvec_0)$ from Eq.~\eqref{eqn:f_U}, compress, and evaluate  $\|f(\Uvec_{0})\|$,  $\|\mathcal{C}\|$ 
    \For{$k = 0,1,\dots, n_{iter}-1$}
        \State \textbf{Step 2:} Compute Jacobian $J_k$ from Eq.~\eqref{eqn:Jacobian} and compress
        
        \State \textbf{Step 3:} Solve linear system Eq.~\eqref{eqn:dU_Newton} for the correction $\Uvec'_k$ using two-site DMRG solver and compress
        \State \quad $\Uvec'_k = \texttt{DMRG}(J_k,-f(\Uvec_k),\alpha,n_{sweeps},\eps_{\mathrm{DMRG}},\chi).compress(\eps_{\mathrm{TT},k},\chi)$
        
        \State \textbf{Step 4:} Update solution via line search  
        \State $\omega_k = 1$
        \For{$m = 1,\dots, n_{line}$}
        \State Compute $\Uvec_{k+1}$ using Eq.~\eqref{eqn:Uvec_kp1} and compress
        \State Evaluate the residue $f(\Uvec_{k+1})$ from Eq.~\eqref{eqn:f_U}, compress, and evaluate $\|f(\Uvec_{k+1})\|$ 
            \If{$\|f(\Uvec_{k+1})\| < \|f(\Uvec_{k})\|$}
                \State Compute $\|\Uvec_{k+1}\|$
                \State \textbf{Break}
            \EndIf
            \State $\omega_k = \omega_k * s$
        \EndFor
        \If{$\|f(\Uvec_{k+1})\| \geq \beta \|f(\Uvec_{k})\|$}
            \State $\eps_{\mathrm{TT},k} = 0.8 \, \eps_{\mathrm{TT},k}$
            \State $\eps_{\mathrm{TT},k} = \max(\eps_{\mathrm{TT},k}, \eps_{\mathrm{TT}})$
        \EndIf
        
        \State \textbf{Check Convergence:}
        \If{$\|f_{k+1}\|/\|\mathcal{C}\| < \eps_{res}$ \textbf{or}
            $\|\Uvec_{k}^{'}\|/\|\Uvec_{k+1}\| < \eps_{cor}$}
            \State \textbf{Break}
        \EndIf
    \EndFor
\end{algorithmic}
\end{algorithm}

\noindent
\begin{algorithm}[H]
\caption{Newton method for solving Eq.~\eqref{eqn:f_U} using Eq. \eqref{eqn:U_Newton}}
\label{algo:Newton_right}
\begin{algorithmic}[1]
    \State \textbf{Input:} Initial guess $\Uvec_0$, TT-rounding inputs $\eps_{\mathrm{TT}}, \chi$ Newton solver inputs $n_{iter}, \eps_{res}, \eps_{cor}$, line search parameters $s, n_{line}$, linear system solver inputs $\eps_{\mathrm{DMRG}}, n_{sweeps}$, Tikhonov regularization parameter $\alpha$
    \State \textbf{Output:} Solution $\Uvec_{k+1}$

    \State Set $\eps_{\mathrm{TT},k}=10^{-3}$ 
    \State \textbf{Step 1:} Compute initial Residue $f(\Uvec_0)$ from Eq.~\eqref{eqn:f_U}, compress, and evaluate  $\|f(\Uvec_{0})\|$,  $\|\mathcal{C}\|$ 
    \For{$k = 0,1,\dots, n_{iter}-1$}
        \State \textbf{Step 2:} Compute Jacobian $J_k$ from Eq.~\eqref{eqn:Jacobian} and compress        
        \State \textbf{Step 3:} Compute and compress $\mathcal{R}_k$ defined in Eq.~\eqref{eqn:U_Newton_line_solve}
        \State \textbf{Step 4:} Solve linear system Eq.~\eqref{eqn:U_Newton_line_solve} for $\widetilde{\Uvec}_{k+1}$ using two-site DMRG solver and compress
        \State \quad $\widetilde{\Uvec}_{k+1} = \texttt{DMRG}(J_k,\mathcal{R}_k),\alpha,n_{sweeps},\eps_{\mathrm{DMRG}},\chi).compress(\eps_{\mathrm{TT},k},\chi)$
        
        \State \textbf{Step 4:} Update solution via line search  
        \State $\omega_k = 1$
        \For{$m = 1,\dots, n_{line}$}
        \State Compute $\Uvec_{k+1}$ using Eq.~\eqref{eqn:U_Newton_line_update} and compress
        \State Evaluate the residue $f(\Uvec_{k+1})$ from Eq.~\eqref{eqn:f_U}, compress, and evaluate $\|f(\Uvec_{k+1})\|$ 
            \If{$\|f(\Uvec_{k+1})\| < \|f(\Uvec_{k})\|$}
                \State Compute $\|\Uvec_{k+1}\|$
                \State \textbf{Break}
            \EndIf
            \State $\omega_k = \omega_k * s$
        \EndFor
        \If{$\|f(\Uvec_{k+1})\| \geq \beta \|f(\Uvec_{k})\|$}
            \State $\eps_{\mathrm{TT},k} = 0.8 \, \eps_{\mathrm{TT},k}$
            \State $\eps_{\mathrm{TT},k} = \max(\eps_{\mathrm{TT},k}, \eps_{\mathrm{TT}})$
        \EndIf
        
        \State \textbf{Check Convergence:}
        \If{$\|f_{k+1}\|/\|\mathcal{C}\| < \eps_{res}$ \textbf{or}
            $\|\Uvec_{k}^{'}\|/\|\Uvec_{k+1}\| < \eps_{cor}$}
            \State \textbf{Break}
        \EndIf
    \EndFor
\end{algorithmic}
\end{algorithm}

\iftrue
\subsection{Compression accuracy}
We introduce an adaptive compression strategy that aligns TT-rounding error with the numerical discretization error by tying it to the mesh spacing. The TT-vectors $\Uvec$, $f(\Uvec)$, and the Jacobian TT-matrix $J$ are compressed with different thresholds. 

Let $\Uvec_\text{exact}$ be the analytical solution of the continuous PDE, $\Uvec_d$ the numerical solution of the discretized PDE, and $\Uvec_{TT}$ its TT representation. Assume $\Uvec_d$ has formal accuracy $\mathcal{O}((\Delta x)^r,(\Delta t)^s)$ and define $C_m \coloneqq \min((\frac{\Delta x}{\Delta x_0})^r, (\frac{\Delta t}{\Delta t_0})^s)$, with $\Delta x_0$ and $\Delta t_0$ stand for the reference mesh spacings, such that the reference mesh yields $C_m = 1$; We have
\begin{align*}
\mathcal{E}\coloneqq \Uvec_\text{exact} - \Uvec_{TT} =  \underbrace{\Uvec_\text{exact} - \Uvec_{d}}_{\mathcal{E}_d} + \underbrace{\Uvec_{d} - \Uvec_{TT}}_{\mathcal{E}_r}
\end{align*}
where, $\mathcal{E}_d \sim \mathcal{O}((\Delta x)^r,(\Delta t)^s)$ and $\mathcal{E}_r$ are the numerical discretization errors, and TT-rounding, respectively. 

From Eq. \eqref{eqn:dU_Newton} we have,
\begin{align*}
    \|\Uvec'_k\| &\leq \|J^{-1}_k\| \, \|f(\Uvec_k)\|, \\
    \|\Uvec'_k\| &\leq \sigma_{\min}^{-1} (J_k) \, \|f(\Uvec_k)\|, \\
    \|\Uvec'_k\| &\lesssim \varepsilon_J^{-1} \, \|f(\Uvec_k)\|.
\end{align*}
Each Newton iteration has the TT-variables compressed as below (the suffix indicates the respective variable):
\begin{align*}
\varepsilon_{TT} &= C_m \, \widetilde{\varepsilon}_{TT} = \mathcal{O}(\Delta t^r, \Delta x^r)\\
\varepsilon_\Uvec &= \varepsilon_{\Uvec'} = \varepsilon_J = \varepsilon_{TT} \\
\varepsilon_f &= \varepsilon_{\mathcal{R}} = \sigma_{\min}(J) \,\varepsilon_{r} \approx \eps_J \eps_r \approx 0.25 \eps_r.
\end{align*}
\fi
\subsection{Multi-level strategy}
We propose a multi-level method which leverages a hierarchy of computational grids—ranging from coarse to fine—to accelerate convergence and improve robustness of the Newton method. A refinement ratio (ratio of mesh spacing of two successive levels) of two is used in both space and time.
%

In this approach, the nonlinear problem is first solved approximately on a coarse grid, where computations are cheaper. The coarse-grid solution is then prolongated (interpolated) to a finer grid 
to serve as an initial guess for the Newton solver at that level. This process continues recursively across multiple levels, from coarse to fine, with each finer-level Newton solve benefiting from a progressively better-informed starting point.

By combining mesh refinement with Newton's method, this strategy mitigates the sensitivity to poor initial guesses, enhances convergence, and can significantly reduce computational cost compared to solving the nonlinear problem directly on the finest grid from scratch.

Consider a hierarchy of $n_l$ space-time meshes, $\Omega^{(l)}$ for $l=0,1,\dots,n_l-1$, with a refinement ratio of two. Each level $l$ has a mesh size of $(\Delta x, \Delta t)2^l$, where $\Delta x$ and $\Delta t$ are the space and time steps at the finest level ($l=0$).
Let the discrete nonlinear system of equations Eq. \eqref{eqn:f_U} derived from $\Omega^{(l)}$ be, 
\begin{align*}
f^{(l)} = f^{(l)}(\Uvec^{(l)}) = 0, \label{eqn:f_U^l}
\end{align*}
where, $f^{(l)}: \mathbb{R}^{N^{(l)}} \rightarrow \mathbb{R}^{N^{(l)}}$ is the nonlinear operator on mesh $\Omega^{(l)}$, $\Uvec^{(l)} \in \mathbb{R}^{N^{(l)}}$ is the space-time solution vector on $\Omega^{(l)}$, $N^{(l)}=(N_x/2^l) \times (N_t/2^l)$ is the number of degrees of freedom on $\Omega^{(l)}$.
Algorithm \ref{algo:Newton_multilevel} describes the steps in the multi-level Newton method.

\begin{algorithm}[H]
\caption{Multilevel space--time solver for Eq.~\eqref{eqn:f_U}}
\label{algo:Newton_multilevel}
\begin{algorithmic}[1]
\State \textbf{Input:} Initial condition $U^0(x)$; number of levels $n_{\ell}$;
prolongation operators $\mathcal{P}^{(\ell,\ell+1)}$;
TT parameters $\eps_{\mathrm{TT}}, \chi$;
Newton parameters $n_{\mathrm{iter}}, \eps_{\mathrm{res}}, \eps_{\mathrm{cor}}, s, n_{\mathrm{line}}$;
DMRG parameters $\eps_{\mathrm{DMRG}}, n_{\mathrm{sweeps}}$;
Tikhonov regularization parameter $\alpha$
\State \textbf{Output:} Space--time solution $\Uvec^{(0)}$ on the finest level

\For{$\ell = n_{\ell}-1, \dots, 0$} 
\Comment{coarsest level $\ell=n_{\ell}-1$ to finest $\ell=0$}
    \State Sample initial condition on $\Omega^{(\ell)}$ and form TT vector $(U^0)^{(\ell)}$
    \State Assemble TT operators $A^{(\ell)}, B^{(\ell)}, \mathcal{C}^{(\ell)}$
    
    \If{$\ell = n_{\ell}-1$} \Comment{coarsest level}
        \State Construct initial guess: $\Uvec^{(\ell)}_0 = \bm{1}_{n_t^{(\ell)}} \otimes (U^0)^{(\ell)}$, where $\bm{1}_{n_t^{(\ell)}} = (1, \dots, 1)^T \in \mathbb{R}^{n_t^{(\ell)}}$
    \Else
        \State Prolongate solution from previous level: $\Uvec^{(\ell)}_0 = \mathcal{P}^{(\ell,\ell+1)} \Uvec^{(\ell+1)}$
    \EndIf
    
    \State Solve the nonlinear problem Eq.~\eqref{eqn:f_U} on level $\ell$ using
    either Algorithm~\ref{algo:Newton_left} or Algorithm~\ref{algo:Newton_right} with initial guess $\Uvec^{(\ell)}_0$
    
    \State Store converged solution $\Uvec^{(\ell)}$
\EndFor

\State \textbf{Return} $\Uvec^{(0)}$
\end{algorithmic}
\end{algorithm}

A prolongation operator $P$ is used to transfer a solution from a coarse grid (level $l$) to a finer grid (level $l-1$), refining the approximation. 
Given a coarse-grid solution $\Uvec^{(l)}$, the prolongation process is represented as $\Uvec^{(l-1)} = P^{(l-1,l)} \Uvec^{(l)}$ where $\Uvec^{(l-1)}$ is the interpolated fine-grid solution. 
A bilinear interpolation is employed in two dimensions wherein fine-grid values are computed using weighted averages of neighboring coarse-grid values (see \cite{Mohr:2004} and pages 68-70 of \cite{Wesseling:1992}). Tensor train representation of the space-time prolongation operator is obtained by the Kronecker product of the 1D linear prolongation operators \cite{kazeev2012low} in space ($P_x$) and time ($P_t$) as below:
\begin{align}
   P^{(l-1,l)} &= P_t^{(l-1,l)} \otimes P_x^{(l-1,l)} = P_t^{(l-1,l)} \bowtie P_x^{(l-1,l)}. \label{eqn:prolonation}
\end{align}

\subsection{Computational complexity}
We assume implicit Euler time discretization and second-order central difference for space derivatives for all the three methods in the following analysis.
Assume spatial size $N_x$ and time size $N_t$, with $q_x=\log_2 N_x$, $q_t=\log_2 N_t$ and let $q=q_x+q_t$. Let $k_N$ be the number of Newton iterations. For QTT, denote typical ranks of vectors/iterates by $r$ and let $\chi=\max(r)$ be the bond dimension defined as the maximum rank of the solution across all cores.

\subsubsection{Classical time stepping}
In classical time stepping approach, we use Newton method in each time step and employ a dense direct solver to solve a linear system for each Newton iteration.
Each Newton iteration solves a dense linear system of $N_x$ unknowns. 
Assuming $k_N$ Newton iterations per time step, cost of $N_t$ time steps is:
\[
\text{Work} \;=\; \mathcal{O}\!\big(k_N\,N_t\,N_x^3\big),
\qquad
\text{Memory} \;=\; \mathcal{O}(N_x^2).
\]

\subsubsection{Single-level space-time QTT}
The space-time solution vector $\Uvec \in\mathbb{R}^{N_x N_t}$ is represented in QTT with $q=q_x+q_t=\log_b(N_x N_t)$ binary modes (each mode with size $b=2$). DMRG method involves solution of a square dense linear system of $b r^2$ unknowns per each site per sweep for which we use the SVD-based least-squares method (scales cubically) with a cost of $\mathcal{O}(b^3 r^6)$. If the ranks $r$ are uniformly bounded by $\chi$ across all DMRG iterations, one DMRG-based linear solve within Newton with \emph{$k_D$ sweeps} costs (up to constants)
\[
\text{Work per Newton step} \;=\; \mathcal{O}\!\big(k_D q\,\chi^6\big),
\qquad
\text{Memory} \;=\; \mathcal{O}\!\big(q\,\chi^2\big).
\]
Total:
\[
\text{Work} \;=\; \mathcal{O}\!\big(k_N\,k_D\,q\,\chi^6\big)  \;=\; \mathcal{O}\!\big(k_N\,k_D\,\log_2(N_x\,N_t)\,\chi^6\big).
\]

\subsubsection{Multi-level space-time QTT}
Assume $k_N$ Newton iterations on each level $\ell=0,1,\ldots,L-1$ with each level have a mesh of size $(2^{q_x-l})\times (2^{q_t-l})$ so that each level has $q^{(\ell)}=q-2\ell$ cores where $q$ is the number of cores on the finest level. Typically, $L=\min(q_t, q_x) - 2$ with $q_t, q_x \geq 2$.
Let $\chi^{(\ell)}$ be the maximum QTT rank on level $\ell$. Each Newton iteration uses $k_D$ DMRG sweeps. Then the total work (up to constants) is
\[
\text{Work} \;=\; 
\sum_{\ell=0}^{L-1} k_N\,k_D\,\mathcal{O}\!\big(q^{(\ell)}\!\,(\chi^{(\ell)})^6\big) = \sum_{\ell=0}^{L-1} k_N\,k_D\,\mathcal{O}\!\big((q-2{\ell})\!\,(\chi^{(\ell)})^6\big)
\]
and the peak memory is dominated by the finest level:
\[
\text{Memory} \;=\; \mathcal{O}\!\big(q\,(\chi^{(0)})^2\big).
\]
If ranks are nondecreasing and uniformly bounded by $\chi$ across levels, and each refinement increases $q^{(\ell)}$ by $2$, then
\begin{align*}
\text{Work} &\;=\; \mathcal{O}\!\big(\,k_N\,k_D\,(q-L+1)L\,\chi^6\big) \;=\; \mathcal{O}\!\big(\,k_N\,k_D\,(\log_2(N_x\,N_t)-L+1)L\,\chi^6\big),
\\
\text{Memory} &\;=\; \mathcal{O}\!\big(q\,\chi^2\big).
\end{align*}

Note that the constant $b$ from binary physical dimension is absorbed into the big-$\mathcal{O}$ terms and the multi-level strategy proceeds strictly from coarse to fine; there is no back-and-forth cycling.
\section{Numerical experiments with nonlinear PDEs}\label{sec:results}
Next, we compare multi-level TT-space-time (ML), single-level TT-space-time (SL), and classical time-stepping (CT) methods on nonlinear initial-boundary value problems, using identical numerical discretization schemes across all three. The numerical error ($e$) and relative error ($e_r$) are quantified through $L^2-$ norm defined as,
\begin{align}
    \|e\|_{L^2}  = \left(\int_\Omega |e_j|^2  ~d\Omega \right)^{1/2} \approx \left(\sum_{j=1}^{N_x} |e_j|^2 \Delta x \right)^{1/2}, \quad \|e_r\| \coloneqq \frac{\|e\|_{L^2}}{\|u\|_{L^2}}
\end{align}
Unless specified otherwise, first-order implicit Euler discretization is used for time derivatives and second-order central discretization is used for space derivatives. Tensor train representation of all the discrete operators used in the following tests is presented in Appendix \ref{sec:tt_operators}.

The solution accuracy, convergence, and runtime depend on multiple parameters, including space and time steps $(\Delta x, \Delta t)$, formal accuracy of the discretization, rounding errors $\eps_{\mathrm{TT}}$ and $\eps_{\mathrm{DMRG}}$, convergence criteria, maximum TT-rank $\chi$, and the
number of mesh levels. 
The truncation errors $\eps_{\mathrm{TT}}$ and $\chi$ should be chosen to be smaller than the discretization error. Often, these quantities are not known a priori and must be estimated empirically. 

In the sequel, we shall present numerical evidence of the convergence of our approach for a variety of types of nonlinear PDEs. It is worth noticing that the need for regularization, in the numerical resolution, arises predominantly in nonlinear hyperbolic problems, where the Jacobian of the monolithic space--time system becomes severely ill-conditioned.
Physically, hyperbolic dynamics are characterized by advective transport, steep gradients, and the formation of sharp features such as shocks or wave fronts.
In a global space--time discretization, these features induce strong anisotropy, leading to Jacobians with rapidly decaying singular spectra. 
Particularly, the ill-conditioning worsens as the system size increases.

Within the tensor-train framework, this ill-conditioning manifests at the level of the local two-site DMRG problems.
The effective local matrix \( A_{\mathrm{loc}} \), obtained by projecting the global Jacobian onto a low-rank TT manifold, may become nearly rank-deficient for stiff nonlinear problems such as the soliton propagation simulated via KdV and sine-Gordon equations. These effects become more prominent as the system size increases.
As a result, small singular values appear in \( A_{\mathrm{loc}} \), rendering the corresponding least-squares problems ill-posed and causing Newton updates to amplify numerical noise rather than reduce the residual.

This behavior is fundamentally different from that observed in parabolic problems, where diffusion provides intrinsic regularization of the operator.
In those cases, the Jacobian remains well-conditioned across mesh refinements, the TT ranks remain moderate, and the local DMRG systems are well-posed without additional stabilization.
For hyperbolic problems, however, the combination of strong nonlinearity and sharp space--time features necessitates explicit regularization to stabilize the local solves and ensure reliable Newton convergence.

All numerical experiments are performed using a Python implementation (Python~3.12.0) on a 2015 MacBook Pro equipped with a 2.5\,GHz Intel Core i7 processor and 16\,GB of 1600\,MHz DDR3 memory.

\subsection{Nonlinear parabolic PDE: Fisher-KPP equation}
The Fisher-Kolmogorov-\\-Petrovsky-Piskunov (Fisher-KPP) equation is a standard example of a reaction-diffusion equation.
Originally developed by Fisher \cite{Fisher:1937} to model the spatial propagation of advantageous alleles in population genetics, this equation has become fundamental across diverse fields, including ecology (species invasion dynamics), combustion theory (flame front propagation), epidemiology (disease spread), and wound healing processes.
\begin{equation}
\frac{\partial u}{\partial t} = D\frac{\partial^2 u}{\partial x^2} + ru(1-u), \quad x \in [-L/2, L/2], \quad t > 0 \label{eqn:fkpp}
\end{equation}
where, $u(x,t)$ represents the population density or concentration, $D>0$ is the diffusion coefficient, and $r>0$ is the intrinsic growth rate. 
The linear and the nonlinear operators are $\Ld = -D \delsux{}{x} -r$, and $\Nd(u) = r u^2$, respectively. The spatial derivatives are discretized using a second-order central difference scheme, while implicit Euler time discretization is used for time discretization. 

The equation admits traveling wave solutions $u(x,t) = U(z)$ with $z = x - ct$ connecting the unstable state $u = 0$ to the stable state $u = 1$, satisfying $\lim_{z \to +\infty} U(z) = 0$ and $\lim_{z \to -\infty} U(z) = 1$. 
While no general closed-form solution exists, Ablowitz and Zeppetella \cite{Ablowitz:1979} found an exact analytical solution for the special case $c = 5/\sqrt{6}$ (with $D=r=1$):
\begin{align}
U(z) = \frac{1}{\left[1 + A e^{z/\sqrt{6}}\right]^2}    \label{eqn:fkpp_analytic}
\end{align}
where $A > 0$ determines the wave position.
We choose $A=1, L=40$, initial condition from Eq. \eqref{eqn:fkpp_analytic} with $t=0$, and boundary conditions $u(-L/2,t)=1, u(L/2,t)=0$ and compute the numerical error w.r.t. the analytical solution.

\begin{figure}[ht]
 \centering
\begin{minipage}{0.5\linewidth}\centering
 \includegraphics[width=2.4in,trim={0.1in 0.in 0.in 0.in},clip]{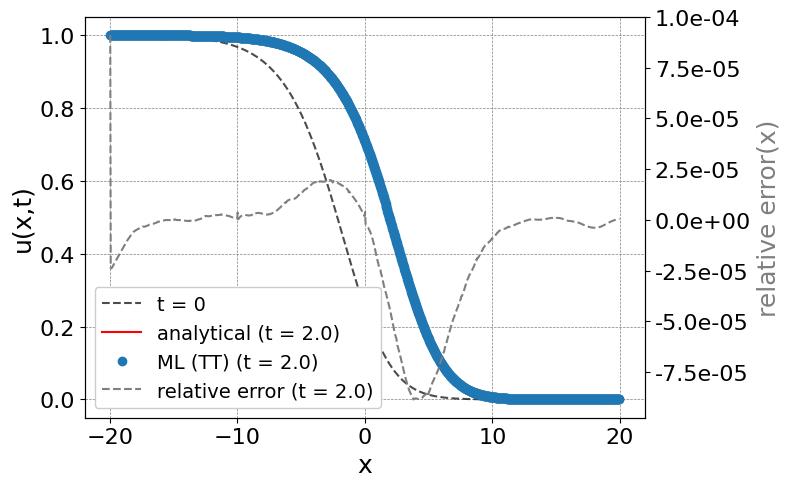}\\(a)
 \end{minipage}
\begin{minipage}{0.45\linewidth}\centering
 \includegraphics[width=2.3in,trim={0.1in 0.in 0.in 0.in},clip]{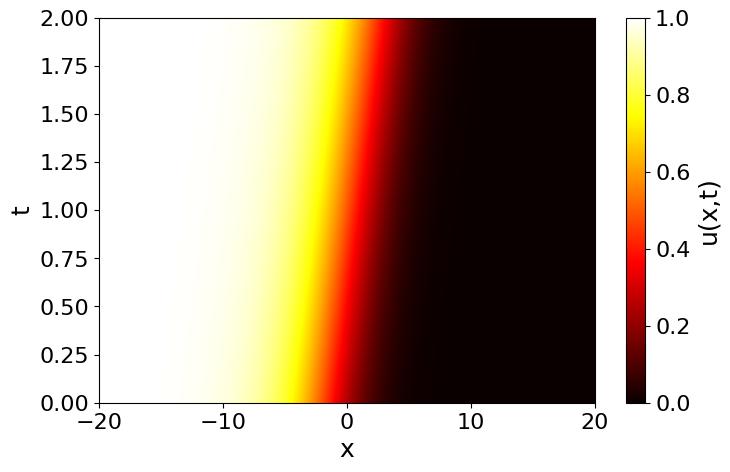}\\(b)
 \end{minipage}
 \\
\begin{minipage}{0.45\linewidth}\centering
 \includegraphics[width=2.3in,trim={0.1in 0.in 0.in 0.in},clip]{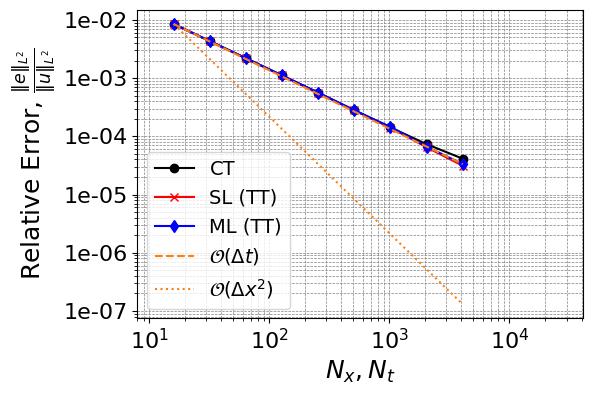}\\(c)
 \end{minipage}
\begin{minipage}{0.5\linewidth}\centering
 \includegraphics[width=2.2in,trim={0.1in 0.in 0.in 0.in},clip]{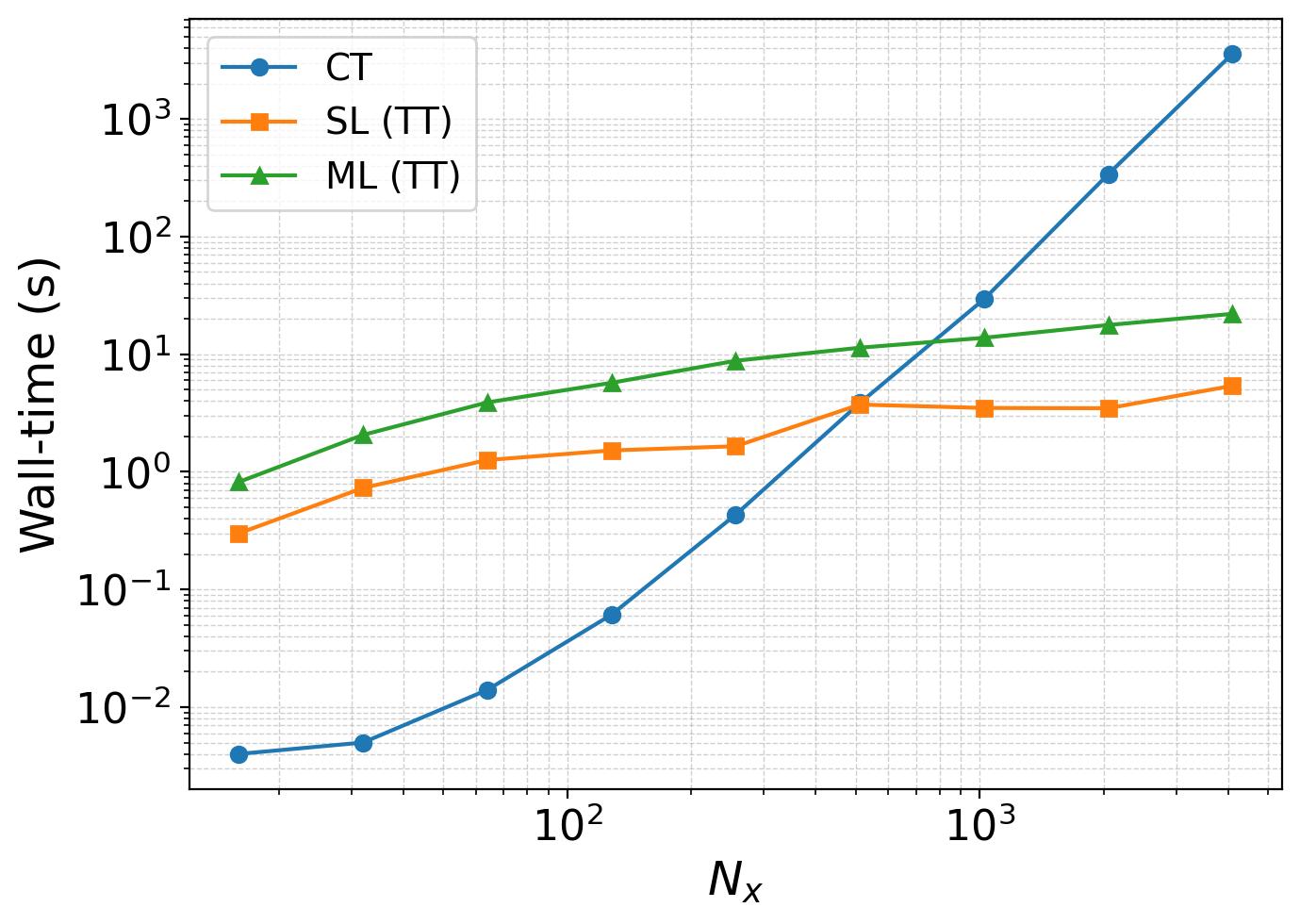}\\(d)
 \end{minipage}
\caption{Traveling wave solution of  Fisher-KPP equation Eq. \eqref{eqn:fkpp}: $D=1, r=1,N_x=N_t=2^{10}, \Delta t=2/N_t$, (a) comparison with analytical solution at $t=2$, , (b) contour of the space-time solution, (c) convergence of numerical error ($\mathcal{O}(\Delta t, \Delta x^2)$) with mesh resolution $N_x, N_t$, and (d) wall-time comparison of classical time-stepping method (CT) and single- (SL) and multi-level (ML) space-time methods.
}
\label{fig:fkpp}
\end{figure}
Fig. \ref{fig:fkpp}(a) shows excellent agreement between the numerical solution with the Multi-Level TT method at $t=2$ on a mesh with $N_x = N_t = 2^{10}$ and $\Delta t = 2/N_t$ and the analytical solution. Fig. \ref{fig:fkpp}(b) presents the solution in the $(x,t)$ plane. Fig. \ref{fig:fkpp}(c) shows that the numerical error at $t=2$ converges as $\mathcal{O}(\Delta t, \Delta x^2)$, consistent with the implicit Euler scheme for time and the second-order central difference for space, which have formal accuracies of $\mathcal{O}(\Delta t)$ and $\mathcal{O}(\Delta x^2)$, respectively.

\begin{table}[]
\centering
\resizebox{\textwidth}{!}{%
\begin{tabular}{@{}c|cr|rrr|cr|ccc|c@{}}
\toprule
$N_x\times N_t$ 
& \multicolumn{2}{c|}{Newton iter} 
& \multicolumn{3}{c|}{wall-time (s)}  
& \multicolumn{2}{c|}{$\max$ rank $\chi$} 
& \multicolumn{3}{c|}{relative error $\|e\|_{L^2}/\|u\|_{L^2}$} 
& $\|e\|_{L^2}$ \\ 
\midrule
& SL & ML 
& SL & ML & CT 
& SL & ML 
& SL & ML & CT 
& CT \\ 
\midrule
$2^4\times 2^4$   
& 3 & 4 
& 0.30 & 0.82 & 0.004 
& 8 & 8 
& 8.51e-03 & 8.51e-03 & 8.51e-03 
& 3.77e-02 \\

$2^5\times 2^5$   
& 4 & 4 
& 0.73 & 2.06 & 0.005 
& 9 & 9 
& 4.32e-03 & 4.32e-03 & 4.32e-03 
& 1.91e-02 \\

$2^6\times 2^6$   
& 4 & 4 
& 1.26 & 3.89 & 0.014 
& 9 & 9 
& 2.19e-03 & 2.19e-03 & 2.19e-03 
& 9.71e-03 \\

$2^7\times 2^7$   
& 4 & 3 
& 1.52 & 5.71 & 0.061 
& 9 & 9 
& 1.11e-03 & 1.11e-03 & 1.11e-03 
& 4.91e-03 \\

$2^8\times 2^8$   
& 4 & 4 
& 1.65 & 8.79 & 0.432 
& 9 & 9
& 5.59e-04 & 5.64e-04 & 5.57e-04 
& 2.47e-03 \\

$2^9\times 2^9$   
& 7 & 4 
& 3.73 & 11.35 & 3.841 
& 9 & 9 
& 2.80e-04 & 2.80e-04 & 2.80e-04 
& 1.24e-03 \\

$2^{10}\times 2^{10}$ 
& 6 & 4 
& 3.49 & 13.78 & 29.329 
& 9 & 9 
& 1.41e-04 & 1.45e-04 & 1.42e-04 
& 6.27e-04 \\

$2^{11}\times 2^{11}$ 
& 6 & 7 
& 3.47 & 17.70 & 338.445 
& 9 & 9 
& 6.47e-05 & 6.63e-05 & 7.33e-05 
& 3.24e-04 \\

$2^{12}\times 2^{12}$ 
& 9 & 8 
& 5.39 & 22.02 & 3570.864 
& 9 & 9 
& 3.12e-05 & 3.34e-05 & 4.12e-05 
& 1.83e-04 \\
\bottomrule
\end{tabular}%
}
\caption{Numerical solution of the Fisher--KPP equation with $D=1$ and $r=1$. 
The TT solver parameters are $\eps_{\mathrm{TT}}=10^{-6}$, $\eps_{\mathrm{DMRG}}=10^{-3}$, and $\eps_{\mathrm{Newton}}=10^{-5}$. The maximum number of DMRG sweeps is fixed to $3$. The multilevel method employs $n_{l}=\min(q_x,q_t)-1$ levels.
Relative errors are defined as $\|e\|_{L^2}/\|u\|_{L^2}$ with $\|u(x,t=2)\|_{L^2}=4.43$. The last column reports the absolute $L^2$ error of the classical time-stepping (CT) solution. SL denotes the single-level TT method, ML the multilevel TT method, and CT classical time stepping.}
\label{tab:fkpp}
\end{table}


Table \ref{tab:fkpp} provides a quantitative comparison of error convergence and computational performance for the three methods: CT, SL, and ML. The results demonstrate that both the Single-Level (SL) and Multi-Level (ML) TT-space-time methods achieve errors nearly identical to the Classical Time-stepping (CT) baseline. This close agreement confirms that the TT-rounding errors (controlled by \(\varepsilon_{TT}\)) are maintained below the level of the spatial and temporal discretization error, ensuring that the low-rank approximation does not compromise solution accuracy.

A key advantage of the Multi-Level approach is its improved nonlinear convergence. As shown in Table \ref{tab:fkpp}, the ML method consistently requires fewer Newton iterations to converge compared to the SL method, particularly on finer grids. This highlights the efficacy of the coarse-grid initialization strategy in providing a high-quality initial guess, thereby enhancing the robustness of the Newton solver.

The computational efficiency of the monolithic TT methods is clearly illustrated in Fig. \ref{fig:fkpp} (d). The wall-time for the CT method exhibits a steep, exponential-like growth with increasing grid resolution, a direct consequence of its \(\mathcal{O}(N_t N_x^3)\) scaling. In stark contrast, both the SL and ML methods, leveraging the logarithmic-scale complexity of the Quantized Tensor Train (QTT) format, show a significantly more favorable log-linear growth in wall-time. This underscores the primary computational benefit of the TT framework: it mitigates the curse of dimensionality inherent in the monolithic space-time approach, making large-scale, high-resolution simulations tractable.

\subsection{Nonlinear parabolic PDE: viscous Burgers' equation}
Consider the nonlinear viscous Burgers' equation below: 
\begin{align}
    u_t + \frac{1}{2}(u^2)_x = \nu u_{xx}, \quad x\in[0,1], \quad t\in[0,1]. \label{eqn:Burgers}
\end{align}
The linear and the nonlinear operators $\Ld = -\nu \delsux{}{x}$, and $\Nd(u) = \frac{1}{2} \delux{}{x}(u^2)$, respectively, are discretized using second-order central difference scheme. 
We consider the following analytical solution \cite{Wood:2006} of Eq. \eqref{eqn:Burgers} to validate the numerical method and assess its accuracy and performance.
\begin{align*}
u(x,t) &= \frac{2\nu \pi \exp{(-\pi^2 \nu t)} \sin(\pi x)}{a + \exp{(-\pi^2 \nu t)} \cos(\pi x)}, \quad 
\quad a > 1, \quad x\in [0,1].
\end{align*}
The initial condition is prescribed from the exact solution at $t=0$, and homogeneous Dirichlet boundary conditions, $u(0,t)=u(1,t)=0$, are imposed for all $t$. The time step is chosen as $\Delta t = 1/N_t$.


We define Reynolds number as $Re = {u_{\max} L}/{\nu}= {2\pi L}/{a}$,
where $L$ denotes the length of the spatial domain and
$u_{\max} = {2\nu\pi}/{a}$ is the maximum advection velocity.
Since $a>1$, the maximum Reynolds number attainable for the initial condition considered above is $Re_{\max} = 2\pi L = 2\pi$,
which is $\mathcal{O}(1)$. This indicates that diffusive effects are significant and that the governing PDE operates predominantly in a parabolic regime.

In contrast, a hyperbolic regime requires $Re \gg 1$, which cannot be achieved with the present choice of initial conditions, as the Reynolds number is directly limited by the viscosity. To investigate the hyperbolic regime, we therefore consider a different class of initial conditions in Sec.~\ref{sec:inviscid_burgers}, which permits large Reynolds numbers and leads to the formation of shock waves.

\begin{figure}[ht]
 \centering
\begin{minipage}{0.5\linewidth}\centering
 \includegraphics[width=2.4in,trim={0.1in 0.in 0.in 0.in},clip]{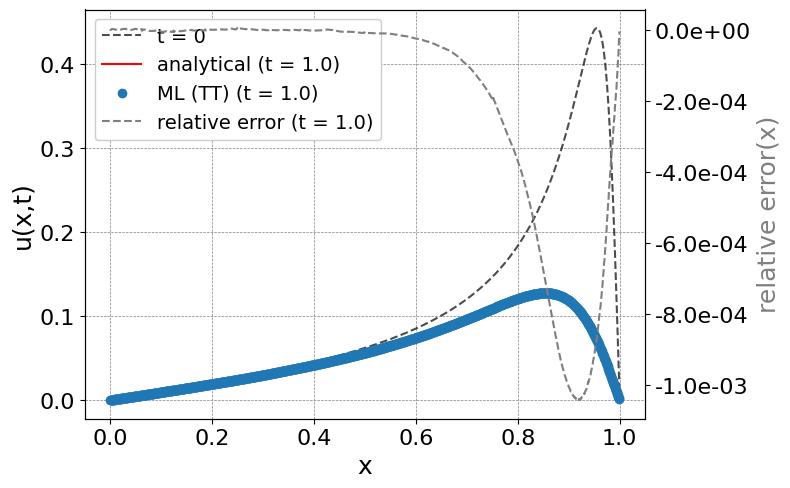}\\(a)
 \end{minipage}
\begin{minipage}{0.45\linewidth}\centering
 \includegraphics[width=2.3in,trim={0.1in 0.in 0.in 0.in},clip]{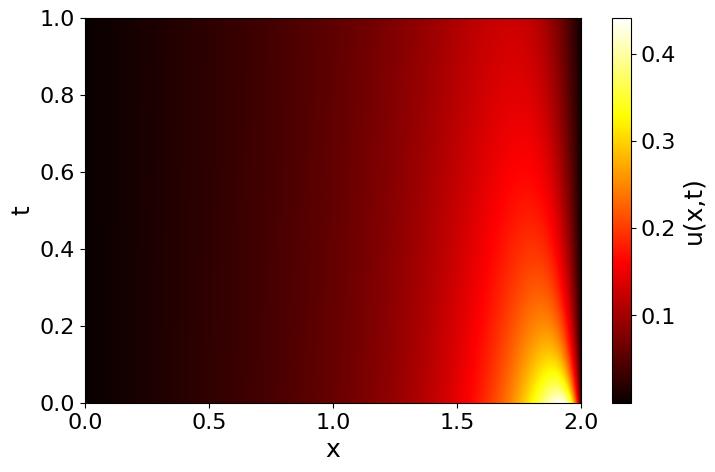}\\(b)
 \end{minipage}\\
\begin{minipage}{0.45\linewidth}\centering
 \includegraphics[width=2.3in,trim={0.1in 0.in 0.in 0.in},clip]{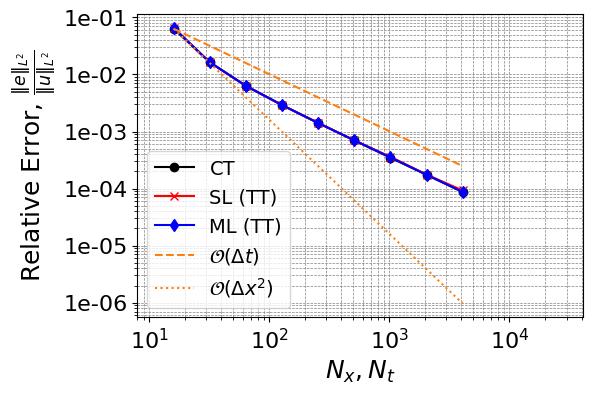}\\(c)
 \end{minipage}
\begin{minipage}{0.45\linewidth}\centering
 \includegraphics[width=2.3in,trim={0.1in 0.in 0.in 0.in},clip]{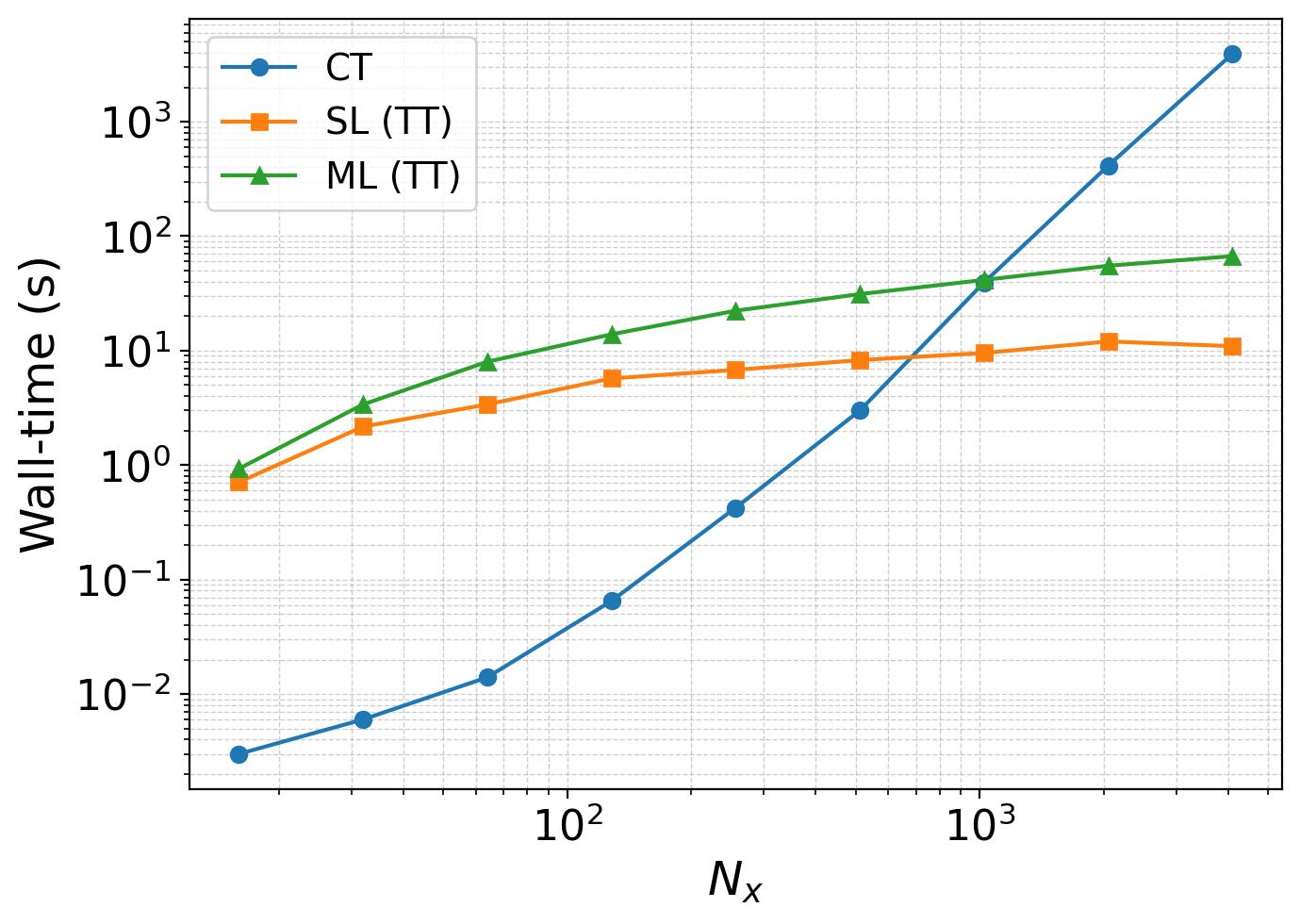}\\(d)
 \end{minipage} 
\caption{Numerical solution of nonlinear viscous Burgers' equation Eq. \eqref{eqn:Burgers} in parabolic regime: $\nu=10^{-2},a=1.01,N_x=N_t=2^{10}, \Delta t=1/N_t$, (a) comparison with analytical solution at $t=1 ~s$, (b) the space-time solution in $xt-$ plane, (c) convergence of numerical error ($\mathcal{O}(\Delta t, \Delta x^2)$) with mesh resolution $N_x, N_t$, and (d) wall-time comparison among classical time-stepping method (CT), single- (SL) and multi-level (ML) space-time methods.}
\label{fig:mps_sin2d}
\end{figure}
Considering $a=1.01$ and $\nu = 0.01$, 
Fig. \ref{fig:mps_sin2d}(a) demonstrates excellent agreement between the numerical ($N_x = N_t = 2^{10}$ and $\Delta t = 1/N_t$) and the analytical solutions at $t=1$. Fig. \ref{fig:mps_sin2d}(b) presents the solution in the $(x,t)$ plane. Fig. \ref{fig:mps_sin2d}(c) shows that the numerical error at $t=1$ converges as $\mathcal{O}(\Delta t, \Delta x^2)$, consistent with the implicit Euler scheme for time and the second-order central difference for space, which have formal accuracies of $\mathcal{O}(\Delta t)$ and $\mathcal{O}(\Delta x^2)$, respectively.

A comparison of computational behavior for the parabolic Burgers equation is summarized in Table \ref{tab:burgers_parabolic}. The TT-based solvers (SL and ML) preserve the discretization accuracy of the classical (CT) method, as their respective errors are virtually identical. This indicates that the tensor compression is lossless with respect to accuracy for the chosen tolerance. The multi-level strategy (ML) provides a clear advantage in nonlinear solver efficiency, consistently converging in fewer iterations. The most significant difference, however, lies in computational scaling. Fig. \ref{fig:mps_sin2d} (d) shows that the wall-time for the CT method surges exponentially with grid refinement, while the SL and ML methods, powered by the compressed QTT representation, exhibit a far more manageable, log-linear increase in runtime.

\begin{table}[]
\centering
\resizebox{\textwidth}{!}{%
\begin{tabular}{@{}c|cr|rrr|cr|ccc|c@{}}
\toprule
$N_x\times N_t$ 
& \multicolumn{2}{c|}{Newton iter} 
& \multicolumn{3}{c|}{wall-time (s)}  
& \multicolumn{2}{c|}{$\max$ rank $\chi$} 
& \multicolumn{3}{c|}{relative error $\|e\|_{L^2}/\|u\|_{L^2}$} 
& $\|e\|_{L^2}$ \\ 
\midrule
$2^{q_x}\times 2^{q_t}$
& SL & ML 
& SL & ML & CT 
& SL & ML 
& SL & ML & CT 
& CT \\ 
\midrule
$2^4\times 2^4$   
& 6 & 4 
& 0.71 & 0.93 & 0.003 
& 7 & 7 
& 6.59e-02 & 6.59e-02 & 6.34e-02 
& 4.65e-03 \\

$2^5\times 2^5$   
& 6 & 5 
& 2.17 & 3.39 & 0.006 
& 9 & 9 
& 1.64e-02 & 1.64e-02 & 1.63e-02 
& 1.16e-03 \\

$2^6\times 2^6$   
& 6 & 5 
& 3.39 & 8.00 & 0.014 
& 10 & 10 
& 6.33e-03 & 6.32e-03 & 6.30e-03 
& 4.47e-04 \\

$2^7\times 2^7$   
& 7 & 4 
& 5.72 & 13.88 & 0.065 
& 10 & 10 
& 2.91e-03 & 2.91e-03 & 2.90e-03 
& 2.05e-04 \\

$2^8\times 2^8$   
& 7 & 5 
& 6.79 & 22.32 & 0.424 
& 10 & 10 
& 1.42e-03 & 1.41e-03 & 1.41e-03 
& 9.98e-05 \\

$2^9\times 2^9$   
& 7 & 4 
& 8.27 & 31.18 & 2.998 
& 10 & 10 
& 7.00e-04 & 7.01e-04 & 7.00e-04 
& 4.94e-05 \\

$2^{10}\times 2^{10}$ 
& 7 & 4 
& 9.54 & 41.47 & 39.368 
& 10 & 10 
& 3.61e-04 & 3.51e-04 & 3.49e-04 
& 2.81e-05 \\

$2^{11}\times 2^{11}$ 
& 8 & 5 
& 12.03 & 55.15 & 412.577 
& 10 & 10 
& 1.76e-04 & 1.75e-04 & 1.74e-04 
& 1.23e-05 \\

$2^{12}\times 2^{12}$ 
& 7 & 5 
& 10.93 & 66.86 & 3912.819 
& 10 & 10 
& 9.26e-05 & 8.71e-05 & 8.70e-05 
& 6.14e-06 \\
\bottomrule
\end{tabular}%
}
\caption{Numerical solution of the Burgers equation with $\nu=0.01$ and $a=1.01$.
The TT solver parameters are $\eps_{\mathrm{TT}}=10^{-6}$, $\eps_{\mathrm{DMRG}}=10^{-3}$, and $\eps_{\mathrm{Newton}}=10^{-5}$.
The maximum number of DMRG sweeps is fixed to $3$, and the multilevel method employs 
$n_{l}=\min(q_x,q_t)-1$ levels.
Relative errors are defined as $\|e\|_{L^2}/\|u\|_{L^2}$ with $\|u(x,t=1)\|_{L^2}=7.06\times10^{-2}$.
The last column reports the absolute $L^2$ error of the classical time-stepping (CT) solution.
SL denotes the single-level TT method, ML the multilevel TT method, and CT classical time stepping.}
\label{tab:burgers_parabolic}
\end{table}

\subsection{Nonlinear hyberbolic PDE: viscous Burgers equation - shock formation} \label{sec:inviscid_burgers}
The viscous Burgers equation, given by Eq. \eqref{eqn:Burgers}, is a fundamental nonlinear partial differential equation that exhibits both convective and diffusive effects, where $\nu > 0$ is the kinematic viscosity. For the initial condition $u_0(x) = -\sin(x)$ on a domain $-\pi\le x\le \pi$, the nonlinear convective term causes the solution to steepen, leading to the formation of shock waves where the solution gradient becomes extremely steep. In the inviscid limit ($\nu \to 0$), the solution develops a true discontinuity, forming a shock that satisfies the Rankine-Hugoniot jump condition. For finite viscosity, the shock is smoothed over a boundary layer of thickness $O(\nu)$, creating a steep but continuous transition. 

\begin{figure}[ht]
 \centering
\begin{minipage}{0.5\linewidth}\centering
 \includegraphics[width=2.2in,trim={0.1in 0.in 0.in 0.in},clip]{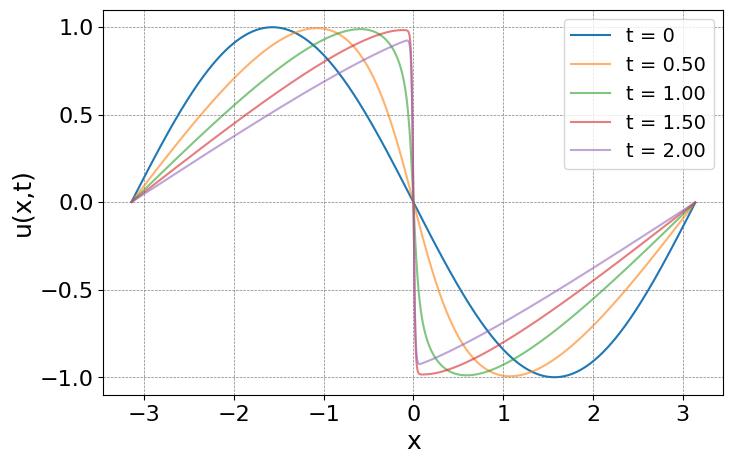}\\(a)
 \end{minipage}
\begin{minipage}{0.45\linewidth}\centering
 \includegraphics[width=2.3in,trim={0.1in 0.in 0.in 0.in},clip]{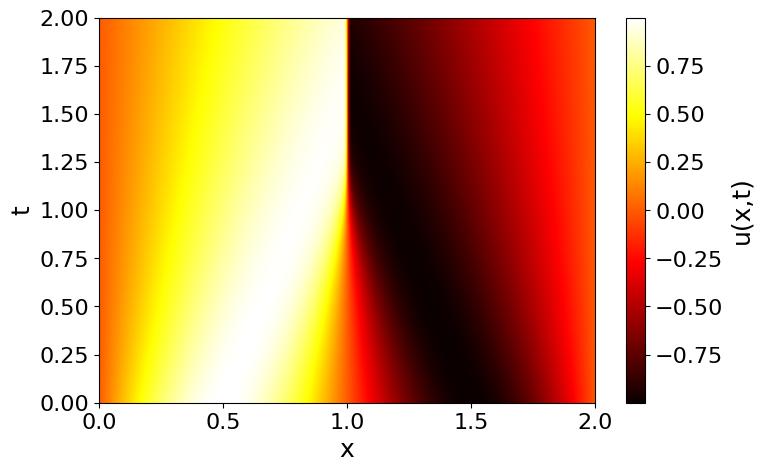}\\(b)
 \end{minipage}
\caption{Evolution of a sine-wave into a shock-wave simulated by solving nonlinear viscous Burgers' equation Eq. \eqref{eqn:Burgers}: $\nu=10^{-2},N_x=N_t=2^{10}, \Delta t=1/N_t$, (a) solutions at every 10 time step interval, (b) contour of the space-time solution. 
}
\label{fig:shock_Burgers}
\end{figure}

Equation \eqref{eqn:Burgers} is solved for $\nu=10^{-2}$ on a domain $x\in[-\pi,\pi]$ discretized with a mesh of size $N_x=N_t=2^{10}$ with $\Delta t=1/N_t$. Fig. \ref{fig:shock_Burgers}(a) shows the evolution of the initial sine wave into a sawtooth-like profile with sharp gradients, demonstrating the competition between nonlinear steepening and viscous diffusion that characterizes the Burgers equation as a simplified model for shock physics in fluid dynamics and traffic flow. Fig. \ref{fig:shock_Burgers}(b) shows the space-time solution wherein the shock is distinguished by a sharp characteristic formed at $t=1$.

\subsection{Nonlinear hyperbolic PDE: sine-Gordon equation}
The sine-Gordon equation is a nonlinear hyperbolic PDE that appears in differential geometry, field theory, and nonlinear wave dynamics, serving as a challenging benchmark for evaluating the accuracy and stability of numerical methods for nonlinear wave equations.
\begin{align}
    &\frac{\partial^2 u}{\partial t^2} - \frac{\partial^2 u}{\partial x^2} + \sin u = 0, 
     \label{eqn:sg}
\end{align}
Equation~\eqref{eqn:sg} is discretized in time using the implicit Euler scheme. The linear operator $\Ld = -\delsux{}{x}$ is discretized using a second-order central finite-difference scheme, and the nonlinear operator
$\Nd(u) = \sin(u)$ is evaluated in TT format via cross approximation
\cite{oseledets2010tt}.
The initial conditions $u(x,0)$ and $\delux{u}{t}(x,0)$ are derived from an analytical kink soliton solution \cite{Bratsos:2008} in the domain $x \in [a, b]$ with $a=-10, b=15$:
\begin{align*}
    &u(x,t) = 4 \arctan\left( \exp\left( \gamma(x - ct - x_0) \right) \right), \\
    &u(a,t) = 0, \quad \delux{u}{x}(b,t) = 0, 
\end{align*}
where \( \gamma = \frac{1}{\sqrt{1 - c^2}} \) is the Lorentz factor with $c$ being the speed of the traveling wave or soliton. 

\begin{figure}[ht]
 \centering
\begin{minipage}{0.5\linewidth}\centering
 \includegraphics[width=2.4in,trim={0.1in 0.in 0.in 0.in},clip]{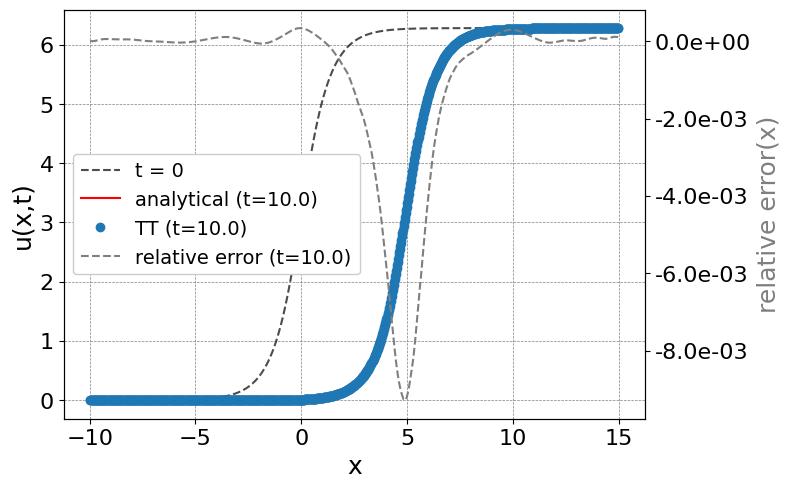}\\(a)
 \end{minipage}
\begin{minipage}{0.45\linewidth}\centering
 \includegraphics[width=2.3in,trim={0.1in 0.in 0.in 0.in},clip]{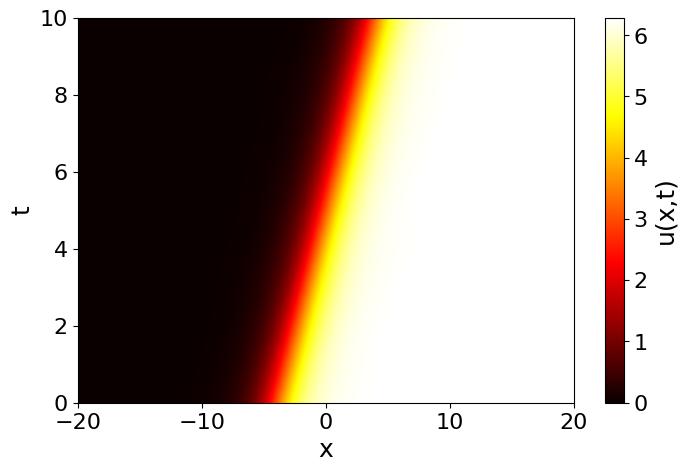}\\(b)
 \end{minipage}
 \\
\begin{minipage}{0.45\linewidth}\centering
 \includegraphics[width=2.5in,trim={0.1in 0.in 0.in 0.in},clip]{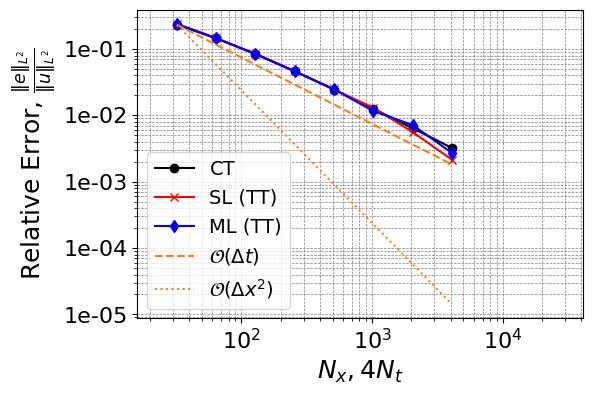}\\(c)
 \end{minipage}
\begin{minipage}{0.5\linewidth}\centering
 \includegraphics[width=2.2in,trim={0.1in 0.in 0.in 0.in},clip]{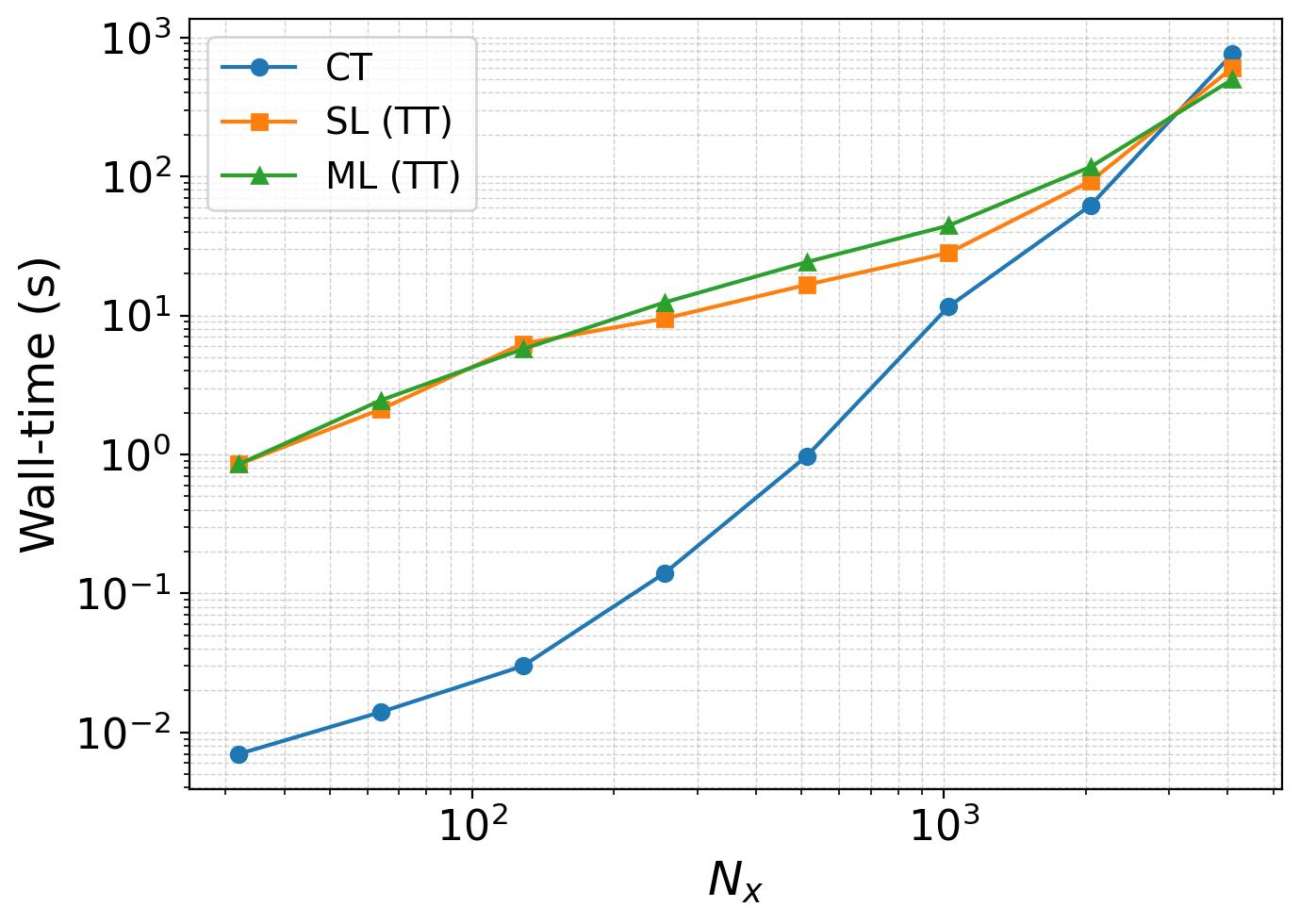}\\(d)
 \end{minipage}
\caption{Numerical solution of the sine--Gordon equation~\eqref{eqn:sg}.
The discretization uses $N_x=2^{10}$ and $N_t=2^{8}$ with time step $\Delta t=10/N_t$.
(a) The numerical and analytical solutions at $t=10$,
(b) space--time contour plot of the numerical solution,
(c) convergence of the numerical error under mesh refinement,
demonstrating first-order temporal accuracy, $\mathcal{O}(\Delta t)$,
due to the use of the implicit Euler time discretization,
(d) wall-clock time of the classical time-stepping method (CT)
with the single-level (SL) and multilevel (ML) space--time methods.}

\label{fig:sg}
\end{figure}
Fig.~\ref{fig:sg} illustrates the numerical solution of the sine--Gordon equation obtained on a
space--time mesh with $N_x = 2^{10}$, $N_t = 2^{8}$, and time step $\Delta t = 10/N_t$, using an
implicit Euler time discretization. The results demonstrate excellent agreement between the
numerical approximation and the analytical solution throughout the simulation interval.

Panel~(a) presents a comparison between the numerical and analytical solutions at the final
time $t=10$. The two profiles are nearly indistinguishable, confirming that the numerical
method accurately captures both the shape and propagation of the kink solution.
Panel~(b) shows a space--time contour plot of the numerical solution, illustrating smooth
temporal evolution and stable soliton propagation without spurious oscillations.

Panel~(c) reports the convergence of the numerical error under mesh refinement. Consistent
with the implicit Euler discretization, the error exhibits first-order convergence with
respect to the time step, $\mathcal{O}(\Delta t)$. For sufficiently fine spatial meshes,
the temporal discretization error dominates, which explains the observed convergence behavior.

Panel~(d) compares the wall-clock time of the classical time-stepping method (CT) with the single-level (SL) and multilevel (ML) space--time solvers. The SL and ML approaches exhibit very similar computational costs across all tested resolutions. This indicates that the additional work associated with the multilevel hierarchy is effectively balanced by the reduction in the number of Newton iterations achieved through coarse-grid initialization, resulting in comparable overall runtimes for SL and ML.
For the largest problem size considered (corresponding to 12 qubits), both SL and ML outperform the classical time-stepping method. This behavior reflects the fundamentally different complexity of the
approaches: the SL and ML space--time solvers exhibit log-linear scaling with respect to the
space--time mesh size, whereas the computational cost of the classical time-stepping method
grows exponentially with problem size.

At the highest mesh resolutions, the computational cost of the TT-based solvers increases
more noticeably. This increase is primarily due to a moderate growth in the TT bond dimension,
which directly affects the cost of the DMRG solver. In particular, the complexity of the DMRG
local solves scales as $\mathcal{O}(\chi^6)$, so even modest increases in the bond dimension
can lead to a measurable rise in runtime. Nevertheless, the TT-based solvers remain competitive
at the resolutions considered in this study. Further reductions in computational cost are
expected from the development of more efficient local solvers and the incorporation of
preconditioning strategies, which could substantially mitigate the dependence on the bond
dimension and extend the favorable scaling behavior to larger problem sizes.

\begin{table}[]
\centering
\resizebox{\textwidth}{!}{%
\begin{tabular}{@{}c|cr|rrr|cr|ccc|c|HHcH@{}}
\toprule
$N_x\times N_t$ & \multicolumn{2}{c|}{Newton iter} & \multicolumn{3}{c|}{wall-time (s)}  & \multicolumn{2}{c|}{$\max$ rank $\chi$} & \multicolumn{3}{c|}{relative error ${\|e\|_{L^2}}/{\|u\|_{L^2}}$} & $\|e\|_{L^2}$ &  \multicolumn{2}{H}{\#time blocks} & $\alpha$ & $\eps_{TT, DMRG}$\\ \midrule
& SL        & ML        & SL            & ML  & CT & SL  & ML & SL  & ML & CT & CT & SL & ML\\ \toprule
$2^{ 5}\times 2^{ 3}$ & 6  & 6 & 0.85 & 0.85 & 0.007 & 6  & 6 & 2.39e-01  & 2.39e-01  &  2.33e-01 & 4.45e+00  & 1 & 1  & $10^{-6}$ & $10^{-3} - 5.12\cdot 10^{-4}$\\
$2^{ 6}\times 2^{ 4}$ & 7  & 4 & 2.11 & 2.45 & 0.014 & 8  & 8 & 1.46e-01  &  1.46e-01 &  1.45e-01 & 2.78e+00  & 1 & 1  & $10^{-6}$ & $10^{-3} - 3.28\cdot 10^{-4}$ \\
$2^{ 7}\times 2^{ 5}$ & 8  & 3 & 6.26 & 5.72 & 0.03 & 9  & 9 & 8.41e-02  & 8.39e-02  &  8.46e-02 & 1.62e+00 & 1 & 1  & $10^{-6}$ & $10^{-3} - 2.62\cdot 10^{-4}$\\
$2^{ 8}\times 2^{ 6}$ & 10 & 4 &  9.49  & 12.41 & 0.14 & 11  & 11 & 4.61e-02  & 4.61e-02  &  4.64e-02 & 8.88e-01 & 1 & 1  & $10^{-7}$ & $10^{-3} - 1.34\cdot 10^{-4}, \chi_{DMRG}=18$\\
$2^{ 9}\times 2^{ 7}$ & 10 & 6 & 16.61  & 24.19 & 0.97 & 13 & 13 & 2.46e-02  & 2.47e-02  &  2.45e-02 & 4.68e-01 & 1 & 1 & $10^{-7}$ & $SL: \text{3e-4}, ML: 10^{-3} - 4.40\cdot 10^{-5}$, $\chi_{DMRG}=18$\\
$2^{10}\times 2^{ 8}$ & 11  & 6  & 28.19 & 44.27 & 11.56 & 17 & 14 & 1.30e-02 & 1.16e-02 &   1.26e-02 & 2.40e-01   & 1 & 1  & $10^{-8}$ & $10^{-3} - 1.15\cdot 10^{-5}$, $\chi_{DMRG}=18$\\
$2^{11}\times 2^{ 9}$ & 30  & 11 & 92.24 & 117.07 & 61.81 & 15 & 15 & 5.61e-03 & 7.15e-03 &  6.41e-03 & 1.22e-01  & 1 & 1 & $10^{-8}$  & $10^{-5}$, $\chi_{DMRG}=18$\\
$2^{12}\times 2^{ 10}$ & 14  & 12 & 603.25 & 498.24 & 761.40 & 17 & 17 & 2.14e-03 & 2.70e-03 &  3.22e-03 & 6.16e-02  & 1 & 1 & $10^{-8}$ & $10^{-6}$, $\chi_{DMRG}=20$
\\
\bottomrule
\end{tabular}%
}
\caption{Numerical results for the sine--Gordon kink--soliton solution.
The TT solver parameters are $\eps_{\mathrm{TT}}=10^{-4}$, $\eps_{\mathrm{DMRG}}=10^{-3}$, and
$\eps_{\mathrm{Newton}}=5\times10^{-4}$. The maximum number of DMRG sweeps is fixed to $2$, and the multilevel method employs $n_{l}=\min(q_x,q_t)-2$ levels.
SL denotes the single-level TT method, ML the multilevel TT method, and CT classical time stepping.
The Tikhonov regularization parameter is denoted by $\alpha$.
Errors are evaluated at $t=10$ and normalized using $\|u\|_{L^2}=19.1$.}
\label{tab:sg}
\end{table}

\subsection{Nonlinear dispersive PDE: Korteweg--de Vries (KdV) equation}
The Korteweg--de Vries (KdV) equation is a fundamental nonlinear partial differential equation that models the evolution of shallow water waves and nonlinear dispersive waves. 
\begin{align}
    &\frac{\partial u}{\partial t} + u \frac{\partial u}{\partial x} + \frac{\partial^3 u}{\partial x^3} = 0, \quad x\in[-L/2, L/2]. \label{eqn:kdv} 
\end{align}
It combines nonlinearity and dispersion in the form of Eq. \eqref{eqn:kdv}, leading to the emergence of solitons—localized wave packets that maintain their shape over time and during interactions. The KdV equation presents notable challenges for numerical simulation due to the stiffness introduced by the third-order derivative and the nonlinear convective term. Capturing soliton dynamics accurately requires careful treatment of both dispersion and nonlinearity, making KdV an excellent 
testbed for validating nonlinear solvers.

Eq. \eqref{eqn:kdv} is discretized with the implicit-Euler time stepping scheme, and the linear operator $\Ld = \delux{}{x^3}$ and the nonlinear operator $\Nd(u) = u\delux{u}{x} = \frac{1}{2}\delux{u^2}{x}$ (recasted into conservative form) are discretized with second-order central differencing schemes.

We consider a spatial domain of length \( L = 30 \) with homogeneous Dirichlet boundary conditions. The initial condition is prescribed by the analytical solution \cite{Holden2012fully,DwivediSarkar:2025} (provided below) for a soliton with speed \( c = 1 \) centered at \( x_0 = -1 \). This solution, which satisfies \( \lim_{x \to \pm\infty} u(x,t) = 0 \), is also used as the reference for quantifying the numerical error.
\begin{align}
    u(x,t) &= 3c\,\mathrm{sech}^2\left( \sqrt{c}\, \frac{(x - ct - x_0)}{2} \right). \label{eqn:kdv_analytical}
\end{align}
\begin{figure}[ht]
 \centering
\begin{minipage}{0.5\linewidth}\centering
 \includegraphics[width=2.4in,trim={0.1in 0.in 0.in 0.in},clip]{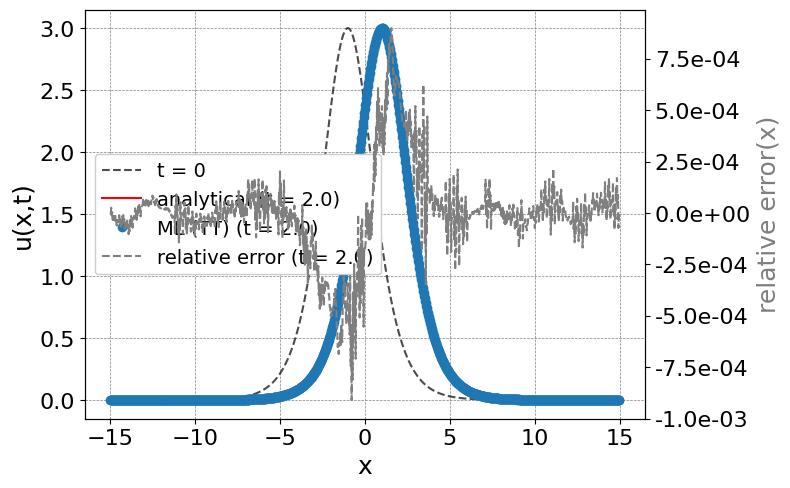}\\(a)
 \end{minipage}
\begin{minipage}{0.45\linewidth}\centering
 \includegraphics[width=2.3in,trim={0.1in 0.in 0.in 0.in},clip]{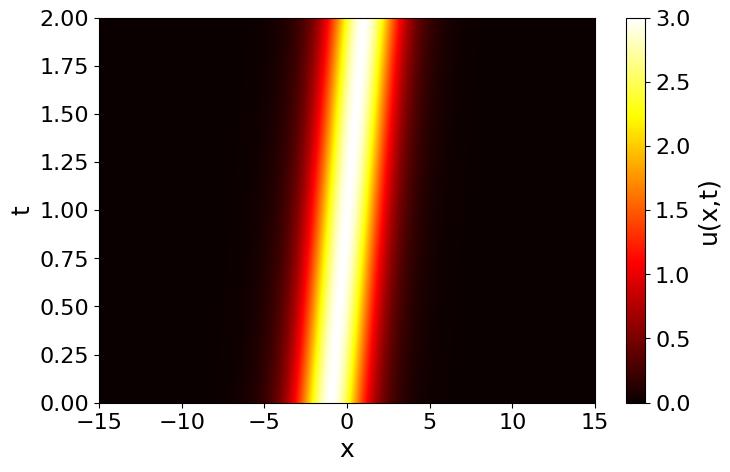}\\(b)
 \end{minipage}\\
\begin{minipage}{0.5\linewidth}\centering
 \includegraphics[width=2.3in,trim={0.1in 0.in 0.in 0.in},clip]{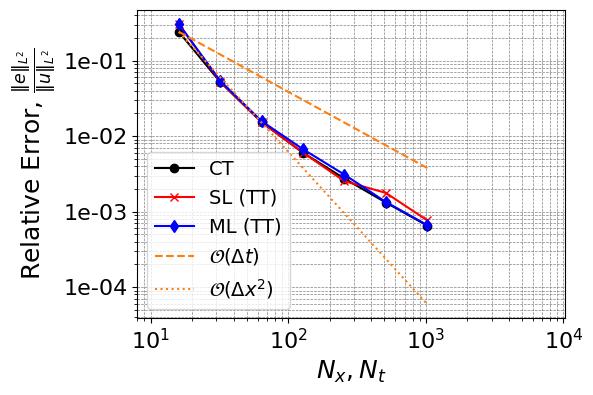}\\(c)
 \end{minipage}
\begin{minipage}{0.45\linewidth}\centering
 \includegraphics[width=2.2in,trim={0.1in 0.in 0.in 0.in},clip]{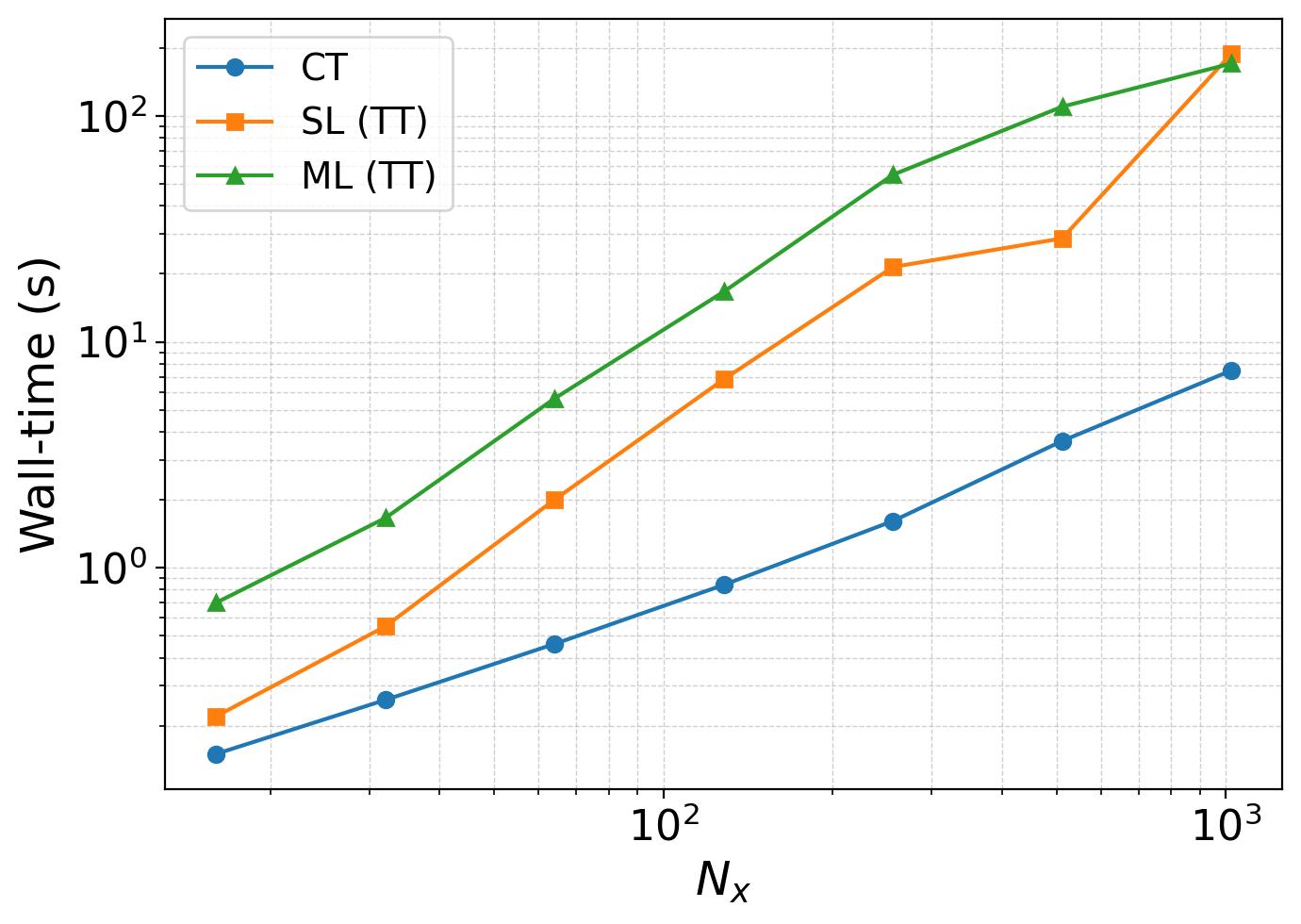}\\(d)
 \end{minipage}
\caption{Numerical solution of the KdV equation~\eqref{eqn:kdv} for the single-soliton case.
The discretization uses $N_x=N_t=2^{10}$ with time step $\Delta t=2/N_t$.
(a) Comparison between the numerical and analytical solutions at $t=2$,
(b) space--time contour plot of the numerical solution,
(c) convergence of the numerical error with mesh refinement:
second-order spatial accuracy $\mathcal{O}(\Delta x^2)$ is observed on coarse meshes,
while for finer meshes the error is dominated by the first-order temporal contribution
$\mathcal{O}(\Delta t)$,
(d) wall-clock time of the classical time-stepping method (CT)
with the single-level (SL) and multilevel (ML) space--time methods.}

\label{fig:kdv_xt}
\end{figure}

Fig.~\ref{fig:kdv_xt} summarizes the numerical results for the KdV equation.
Panel~(a) demonstrates excellent agreement between the numerical solution (ML-TT method) and the analytical
single-soliton solution at $t=2$, obtained using a uniform space--time discretization with
$N_x = N_t = 2^{10}$ and time step $\Delta t = 2/N_t$.
The numerical and analytical profiles are nearly indistinguishable, confirming the accuracy
of the proposed method.

Panel~(b) presents the numerical solution in the $(x,t)$ plane, illustrating the smooth
space--time evolution and stable propagation of the soliton without visible numerical
artifacts.

Panel~(c) shows the convergence behavior of the numerical error evaluated at $t=2$.
The results indicate that the error converges as
$\mathcal{O}(\Delta t, \Delta x^2)$ under mesh refinement.
This behavior is consistent with the use of an implicit Euler scheme for time discretization,
which is first-order accurate, and a second-order central finite-difference scheme
for spatial discretization.

\iftrue
\begin{table}[htbp]
\centering
\resizebox{\textwidth}{!}{%
\begin{tabular}{@{}c|cr|rrr|cr|ccc|c@{}}
\toprule
$N_x\times N_t$ & \multicolumn{2}{c|}{Newton iter} & \multicolumn{3}{c|}{wall-time (s)}  & \multicolumn{2}{c|}{$\max$ rank $\chi$} & \multicolumn{3}{c|}{relative error ${\|e\|_{L^2}}/{\|u\|_{L^2}}$} & $\|e\|_{L^2}$\\ \midrule
& SL        & ML        & SL            & ML  & CT & SL  & ML & SL  & ML & CT & CT\\ \midrule
$2^{ 4}\times 2^{ 4}$ & 3  & 4 & 0.22   & 0.70  & 0.15 & 8  & 7 & 3.07e-01  & 3.07e-01  &  2.42e-01 & 1.18e+00  \\
$2^{ 5}\times 2^{ 5}$ & 3  & 3 &  0.55  & 1.66  & 0.26 & 10  & 10 & 5.29e-02  & 5.42e-02  &  5.23e-02 & 2.56e-01  \\
$2^{ 6}\times 2^{ 6}$ & 3  & 3  & 2.00  & 5.61  & 0.46 & 13 & 13 & 1.53e-02 & 1.58e-02  &   1.53e-02 & 7.50e-02  \\
$2^{ 7}\times 2^{ 7}$ & 3  & 2 &  6.82  & 16.66  & 0.84  & 13  & 13 & 6.05e-03  & 6.71e-03  &  6.03e-03 & 2.95e-02  \\
$2^{ 8}\times 2^{ 8}$ & 5  & 3 & 21.44   & 54.87  & 1.61  & 13  & 13 &  2.53e-03 & 3.09e-03  & 2.74e-03 & 1.34e-02  \\
$2^{ 9}\times 2^{ 9}$ & 5  & 3 &  28.70  & 109.83 & 3.64 & 13  & 13 &  1.77e-03 & 1.34e-03  & 1.32e-03 & 6.49e-03  \\
$2^{10}\times 2^{10}$ & 13  & 3 &  188.14  &  170.50 & 7.47 & 13  & 13 &  7.65e-04 & 6.62e-04  &  6.53e-04 & 3.20e-03  \\
\bottomrule
\end{tabular}%
}
\caption{Numerical results for the KdV single-soliton problem.
The TT solver parameters are $\eps_{\mathrm{TT}}=10^{-6}$, $\eps_{\mathrm{DMRG}}=10^{-3}$, and $\eps_{\mathrm{Newton}}=10^{-3}$. The maximum number of DMRG sweeps is fixed to $3$, and the Newton iteration is capped at $20$ iterations with line-search parameter $s=0.8$. The multilevel method employs $n_{l}=\min(q_x,q_t)-2$ levels. Here, SL denotes the single-level TT method and ML the multilevel TT method. Errors are normalized using $\|u\|_{L^2}=11.2$. The computational domain length is $L=30$ and the soliton speed is $c=1$.
The Tikhonov regularization parameter is set to $\alpha=10^{-12}$.}

\label{tab:kdv_single}
\end{table}
The comparative analysis in Table \ref{tab:kdv_single} confirms two key advantages of the TT-space-time methods for the KdV test case. First, the solution accuracy of the SL and ML solvers comparable (nearly identical)  with the CT reference, demonstrating controlled rounding errors. Second, the ML solver achieves this accuracy with fewer Newton iterations than its SL counterpart, underscoring the benefit of the multi-level initialization. 
\fi

For hyperbolic and dispersive problems, such as the KdV equation, solving the monolithic space–time system in TT-format becomes increasingly challenging as the space–time resolution is refined. In particular, the local linear systems arising in the DMRG-based solution of the Newton correction are often severely ill-conditioned, with conditioning deteriorating rapidly as the space–time mesh is refined. This ill-conditioning leads to stagnation or failure of the Newton iteration unless appropriate stabilization is employed.

To mitigate this difficulty, Tikhonov regularization is applied at the level of the local DMRG solves, stabilizing the projected linear systems while preserving the global low-rank structure. Although this regularization improves robustness in both single-level and multi-level formulations, it is not sufficient on its own to guarantee convergence on fine meshes.

Robust convergence is achieved through the multilevel strategy, in which coarse-grid space–time solutions are used to initialize Newton iterations on successively refined grids. This nonlinear continuation in resolution substantially alleviates the effects of ill-conditioning and yields near mesh-independent convergence. In contrast, single-level space–time TT solvers exhibit rapidly degrading convergence behavior and fail for sufficiently fine discretizations.

Numerically, the multilevel formulation requires approximately the same number of Newton iterations (typically three) across all tested resolutions. By contrast, the single-level approach exhibits a steady increase in Newton iterations with mesh refinement. It fails to converge for larger problem sizes (e.g., space–time meshes exceeding $2^{10}\times 2^{10}$ degrees of freedom).

For all KdV experiments, the TT bond dimension is fixed at 13. While adaptive rank strategies are possible, our results indicate that increased bond dimensions do not lead to improved accuracy for this class of problems. On the contrary, constraining the bond dimension enhances both convergence robustness and computational efficiency. The chosen bond dimension, therefore, provides an effective balance between accuracy and stability.

Provided that TT ranks remain bounded, the computational complexity of the TT-based method remains controlled. In particular, classical solvers exhibit rapid (typically exponential) growth in computational cost under space–time refinement, whereas the TT-based approach displays logarithmic–linear scaling with respect to the total number of space–time degrees of freedom, assuming near-constant Newton iteration counts. This behavior is consistently observed in our experiments: increasing the problem size from $9$ to $10$ qubits results in only a modest increase in runtime, whereas the single-level formulation already shows near-doubling of computational cost at earlier refinements.

It is worth mentioning that, although the present study is restricted to one-dimensional spatial problems, the advantages of the multilevel TT-space–time framework are expected to become more pronounced in higher dimensions. While classical methods typically suffer exponential complexity growth with spatial dimension, TT-based methods can mitigate this effect provided that tensor ranks remain moderate. Extensions to two- and three-dimensional problems are therefore a natural direction for future work.
\section{Conclusion}\label{sec:conclusion}
This work revisits monolithic space–time solvers for nonlinear PDEs from a pragmatic standpoint. While tensor-train representations offer substantial compression of global space–time discretizations, our results make clear that compression alone does not resolve the central difficulty: robust nonlinear convergence. In advection-dominated and hyperbolic regimes, single-level monolithic Newton iterations frequently stagnate, even when low-rank representations are available.

The key contribution of this paper is the introduction of a nonlinear multilevel strategy, embedded directly in the TT framework, whose sole purpose is to construct reliable initial guesses for monolithic Newton solves. Each level is solved independently until convergence and then prolongated to the next refinement. No residual correction or multilevel cycling is performed. This distinction is essential: the method is not designed only to accelerate linear convergence, but also to prevent nonlinear failure. A second essential component is the regularization of the Newton–DMRG iteration. For hyperbolic problems, the projected local systems arising in TT solvers are often nearly singular. Without regularization, Newton updates become unstable and amplify noise. With mild Tikhonov stabilization, convergence is restored while preserving low-rank structure. In this setting, regularization is not optional. All numerical experiments presented here are restricted to one spatial dimension, where classical time stepping already performs well. Even in this favorable regime, the multilevel TT solver exhibits improved robustness and competitive scaling. More importantly, the underlying motivation is dimensional: classical methods scale exponentially with spatial dimension, whereas TT-based formulations retain exponential compression under moderate rank growth. The multilevel strategy proposed here is dimension-agnostic and is expected to play an increasingly important role in higher-dimensional settings.

The effectiveness of monolithic space–time solvers for nonlinear PDEs hinges not on compression alone, but on how Newton’s method is stabilized and initialized in TT framework. The combination of a nonlinear multilevel strategy with a regularized Newton–DMRG iteration provides a viable path forward, particularly for challenging hyperbolic dynamics.

\section*{Acknowledgment}
The last three authors gratefully acknowledge support from grant number KUEXT-8434000491, a collaborative effort between Khalifa University, Technology Innovation Institute, and the Emirates Nuclear Energy Company.
This research was carried out by the Emirates Nuclear Technology Center, a collaboration among Khalifa University of Science and Technology, Emirates Nuclear Energy Company, and the Federal Authority for Nuclear Regulation. 

\appendix
\section{Tensor train representation of the operators}\label{sec:tt_operators}
Let $\mathcal{T}^{(n)}{(l,d,u;a_1,a_2)}$ be a $2^n\times 2^n$ Toeplitz matrix with sub-diagonal $l$, diagonal $d$, and super-diagonal $u$, modified by boundary conditions $a_1, a_2$ (Eq. \eqref{def:aLbL}). For periodic boundaries, replace $\mathcal{T}^{(n)}{(l,d,u;a_1,a_2)}$ with the circulant tridiagonal matrix $\mathcal{T}_c^{(n)}{(l,d,u)}$ of size $2^n\times 2^n$, defined as:
\begin{equation*}
  \mathcal{T}^{(n)}{(l,d,u;a_1,a_2)}= \underbrace{\toeb{l}{d}{u}{a_1l+d}{a_2u+d}}_{2^n \times 2^n} , ~
  \mathcal{T}_c^{(n)}{(l,d,u)}= \underbrace{\toec{l}{d}{u}}_{2^n \times 2^n} .
  \label{def:A_Toeplitz}
\end{equation*}

\noindent Following the procedure in~\cite{kazeev2012low}, the TT representations of $\mathcal{T}^{(n)}(l,d,u;a_1,a_2)$ and $\mathcal{T}_c^{(n)}(l,d,u)$ are:
\begin{align*}
        \mathcal{T}^{(n)}{(l,d,u;a_1,a_2)} &=\begin{bmatrix} I_1 & I & \Jt & J & I_2 \end{bmatrix} \bowtie \begin{bmatrix} I_1 &  &  & & \\ & I & \Jt & J & \\ & & J & & \\ & & & \Jt & \\ & & & & I_2 \end{bmatrix}^{\bowtie (n-2)} \hspace{-0.2in}\bowtie  \begin{bmatrix} a_1 I_1\\d I + u J + l \Jt \\ l J  \\ u \Jt \\ a_2 I_2\end{bmatrix},\\
    \mathcal{T}_c^{(n)}{(l,d,u)} &=\begin{bmatrix} I & \Jt+J & J+\Jt \end{bmatrix} \bowtie \begin{bmatrix} I & \Jt & J \\ & J & \\ & & \Jt \end{bmatrix}^{\bowtie (n-2)} \bowtie  \begin{bmatrix} d I + u J + l \Jt \\ l J  \\ u \Jt \end{bmatrix}.
\end{align*}
where,
\begin{align*}
  {I}&=\begin{pmatrix} 1 & 0 \\  0 & 1 \end{pmatrix},
  {J}=\begin{pmatrix} 0 & 1 \\  0 & 0 \end{pmatrix},
{\Jt}=\begin{pmatrix} 0 & 0 \\  1 & 0 \end{pmatrix},
{I_1}=\begin{pmatrix} 1 & 0 \\  0 & 0 \end{pmatrix},
{I_2}=\begin{pmatrix} 0 & 0 \\  0 & 1 \end{pmatrix}.
  \label{def:QTT_blocks}
\end{align*}
We use $\bowtie$ to denote the rank-core product between two tensors, defined as matrix multiplication with element-wise Kronecker product \cite{kazeev2012low}.
Matrix and TT representations of a pentadiagonal matrix $\mathcal{T}_5^{(n)}(l_2,l_1,d,u_1,u_2;a_{1,1},a_{1,2},a_{1,3},a_{2,1},a_{2,2},a_{2,3}),$ and a circulant pentadiagonal matrix $\mathcal{T}_{c5}^{(n)}(l_2,l_1,d,u_1,u_2)$ of size $2^n\times 2^n$ are:
\begin{align}    
\mathcal{T}^{(n)}_5(l_2,l_1,d,u_1,u_2&;a_{1,1},a_{1,2},a_{1,3},a_{2,1},a_{2,2},a_{2,3}) = \begin{bmatrix*}[r]
a_{1,1}+d   & a_{1,2}+u_1 & u_2 &    &    &  &     \\
a_{1,3}+l_1 & d   & u_1 & u_2 &    &  &      \\
l_2 & l_1 & d   & u_1 & u_2 &  &      \\
 & \ddots & \ddots & \ddots & \ddots & \ddots &  \\
   &  & l_2 & l_1 & d   & u_1    & u_2   \\
   &  &  & l_2 & l_1 & d     &  a_{2,3}+u_1   \\
   &  &  &    & l_2 &  a_{2,2}+l_1   & a_{2,1}+d     \\
\end{bmatrix*}, \\
\mathcal{T}^{(n)}_{c5} &= \begin{bmatrix*}[l]
d   & u_1 & u_2 &    &    & l_2 & l_1    \\
l_1 & d   & u_1 & u_2 &    &  & l_2     \\
l_2 & l_1 & d   & u_1 & u_2 &  &      \\
 & \ddots & \ddots & \ddots & \ddots & \ddots &  \\
   &    & l_2 & l_1 & d   & u_1    & u_2   \\
u_2   &  &  & l_2 & l_1 & d     & u_1   \\
u_1   & u_2 &  &    & l_2 & l_1   & d     \\
\end{bmatrix*},\notag\\
\mathcal{T}^{(n)}_5 &= 
\begin{bmatrix} I_1 & I & \Jt & J & I_2 \end{bmatrix} \bowtie \begin{bmatrix} I_1 &  &  & & \\ & I & \Jt & J & \\ & & J & & \\ & & & \Jt & \\ & & & & I_2 \end{bmatrix}^{\bowtie (n-2)} \hspace{-0.2in} \bowtie  \begin{bmatrix}  a_{1,1} I_1 + a_{1,2} J + a_{1,3} \Jt\\l_1 \Jt + d I + u_1 J \\ l_1 J + l_2 I \\ u_1 \Jt + u_2 I \\ a_{2,1} I_2 + a_{2,2} \Jt + a_{2,3} J\end{bmatrix} \label{def:QTT_pentadiag_bc},\\
\mathcal{T}^{(n)}_{c5}  &= \begin{bmatrix} I & \Jt+J & J+\Jt \end{bmatrix} \bowtie \begin{bmatrix} I & \Jt & J \\ & J & \\ & & \Jt \end{bmatrix}^{\bowtie (n-2)} \bowtie  \begin{bmatrix} l_1 \Jt + d I + u_1 J \\ l_1 J + l_2 I \\ u_1 \Jt + u_2 I\end{bmatrix} \label{def:QTT_pentadiag_periodic},
\end{align}

Matrix representations of the operators $D_t, J_t, K_t, D_x, D_{xx}, D_{xxx}$ are given below:
\begin{align*}
D_t &= \begin{bmatrix*}[r]
1      &        &        &        \\
-1     & 1      &        &        \\
       & \ddots & \ddots &        \\
       &        & -1     & 1
\end{bmatrix*}, 
\quad J_t = \frac{1}{2}\begin{bmatrix*}[r]
1      &        &        &        \\
1     & 1      &        &        \\
       & \ddots & \ddots &        \\
       &        & 1     & 1
\end{bmatrix*}, 
\quad K_t = \frac{1}{4}\begin{bmatrix*}[r]
1      &        &        &    &    \\
2     & 1      &        &     &   \\
1     & 2      &   1     &     &   \\
       & \ddots & \ddots &  \ddots &     \\
       & &  1      & 2     & 1
\end{bmatrix*}, \\
D_x &= \frac{1}{2 \Delta x}\begin{bmatrix*}[r]
-a_1      & 1       &        &        &        \\
-1      & 0     & 1      &        &        \\
       &  \ddots       & \ddots & \ddots &        \\
       &        &  -1     & a_2
\end{bmatrix*}, 
\quad D_{xx} = \frac{1}{\Delta x^2}\begin{bmatrix*}[r]
a_1-2  & 1      &        &   \\
1      & -2     & 1      &        \\
       & \ddots & \ddots & \ddots \\
       &        &  1     & a_2-2
\end{bmatrix*},\\
D_{tt} &= \begin{bmatrix*}[r]
1      &        &        &    &    \\
-2     & 1      &        &     &   \\
1     & -2      &   1     &     &   \\
       & \ddots & \ddots &  \ddots &     \\
       & &  1      & -2     & 1
\end{bmatrix*}, \quad 
D_{xxx}=\frac{d^3 u}{dx^3} \approx \frac{1}{2\Delta x^3}
\begin{bmatrix*}[r]
0 &-2 &  1 &  &  &  &  \\
2 & 0 & -2 & 1 &  &  &  \\
-1 & 2 & 0 & -2 & 1 &  &  \\
 & \ddots & \ddots & \ddots & \ddots & \ddots &  \\
 &  & -1 & 2 & 0 & -2 & 1 \\
 &  &  & -1 & 2 & 0 & -2 \\
 &  &  &  & -1 & 2 & 0 \\
\end{bmatrix*},
\end{align*}
where $a_1$ and $a_2$ denote the boundary condition types at the left and right boundaries, respectively, as defined in Eq. \eqref{def:aLbL}.
TT representations of the operators $D_t, J_t, K_t, D_{tt}, D_x, D_{xx}, D_{xxx}$ are: 
\begin{align*}
    D_t &= \mathcal{T}^{(q_t)}{(-1,1,0;0,0)},\quad
    J_t = \frac{1}{2}\mathcal{T}^{(q_t)}{(1,1,0;0,0)},\quad
    K_{t} = \frac{1}{4}\mathcal{T}_5^{(q_t)}{(1,2,1,0,0)},\\
    D_{tt} &= \mathcal{T}_5^{(q_t)}{(1,-2,1,0,0)},\quad
    D_x =\frac{1}{2\Delta x}\mathcal{T}^{(q_x)}{(-1,0,1;a_1,a_2)}, 
    \\D_{xx} &=\frac{1}{\Delta x^2}\mathcal{T}^{(q_x)}{(1,-2,1;a_1,a_2)},\quad
    D_{xxx} =\frac{1}{2\Delta x^3}\mathcal{T}_{c5}^{(q_x)}{(-1,2,0,-2,1)}.
\end{align*}
\bibliographystyle{siamplain}
\bibliography{tt}

\end{document}

%% file: Misc/shared_siam.tex
\headers{Multi-level TT-space-time solver}{Rapaka \textit{et al.}}

\title{
A Quantum-Inspired Multi-Level Tensor-Train Monolithic Space-Time Method for Nonlinear PDEs
}

\author{
  N. R. Rapaka \footnotemark[5]\,\,\,\thanks{
    Department of Mathematics, College of Computing and Mathematical Sciences,\\
    Khalifa University of Science and Technology, Abu Dhabi, UAE. email: \email{narsimha.rapaka@ku.ac.ae}
  }
  \and
  R. Peddinti\thanks{
    Technology Innovation Institute, Abu Dhabi, UAE\\
    emails: \email{raghavendra.peddinti@tii.ae},
        \email{egor.tiunov@tii.ae},
        \email{leandro.aolita@tii.ae}.
  }
  \and
  E. Tiunov\footnotemark[2]
  \and
  N. J. Faraj\thanks{
 Emirates Nuclear Energy Company,
  Masdar City, Abu Dhabi, UAE. email:\email{noor.faraj@enec.ae}
  }
  \and
  A. N. Alkhooori\footnotemark[3]
  \and
  L. Aolita\footnotemark[2]
  \and
  Y. Addad \footnotemark[1]\,\,\,\thanks{
    Department of Mechanical and Nuclear Engineering,\\
    Khalifa University of Science and Technology, Abu Dhabi, UAE\\
    email: \email{yacine.addad@ku.ac.ae}
  }
  \and
  M. K. Riahi\footnotemark[1]\,
  \thanks{
    Emirates Nuclear Technology Center,\\
    Khalifa University of Science and Technology, Abu Dhabi, UAE\\
    Correspondence should be sent to: \email{mohamed.riahi@ku.ac.ae}
  }
}

\newcommand{\citep}[1]{\cite{#1}}
\input{Misc/packages}
\input{Misc/newcommands}


%% file: Misc/packages.tex
\usepackage{algorithm}
\usepackage{algorithmicx}
\usepackage{algpseudocode} 
\usepackage{amsmath}
\usepackage{amssymb}
\usepackage{slantsc}

\usepackage{diagbox}
\usepackage{xcolor,colortbl}
\usepackage{multirow}
\usepackage{resizegather}
\usepackage{tabu}
\usepackage{epstopdf}
\usepackage{longtable}
\usepackage{caption}
\usepackage{lscape}
\usepackage{graphicx}
\usepackage{latexsym}
\usepackage{wasysym}
\usepackage{fontenc}
\usepackage{rotating}
\usepackage{array}
\usepackage{tikz}

\usepackage{spalign}

\usepackage{mathtools}
\usepackage{booktabs}

%% file: Misc/newcommands.tex
\newcommand{\bm}{\mathbf}

\newsavebox{\astrutbox}
\sbox{\astrutbox}{\rule[-5pt]{0pt}{20pt}}

\newcommand{\be}{\begin{equation}}
\newcommand{\ee}{\end{equation}}
\newcommand{\bee}{\begin{eqnarray}}
\newcommand{\eee}{\end{eqnarray}}

\newcommand {\delux}[2]{\frac{\partial{#1}}{\partial{#2}}}

\newcommand {\delsux}[2]{\frac{\partial^2{#1}}{\partial{{#2}^2}}}

\newcommand {\eps}{\varepsilon}

\definecolor{Gray}{gray}{0.85}
\definecolor{LightCyan}{rgb}{0.88,1,1}
\definecolor{Red}{rgb}{1,.5,0}
\newcolumntype{g}{>{\columncolor{Gray}}c}
\newcolumntype{y}{>{\columncolor{LightCyan}}c}
\newcolumntype{o}{>{\columncolor{Red}}c}

\newcommand{\Jt}{J^{'}}

\newcommand{\toec}[3]{
  {\begin{pmatrix*}[r]
    {#2} & {#3}  &  &  & {#1} \\
    {#1} & {#2} & {#3} &  & \\
     & \ddots & \ddots & \ddots &  \\
      &  &  {#1} & {#2} &  {#3}\\
    {#3} &   &  & {#1} &  {#2} 
  \end{pmatrix*}}
}
\newcommand{\toeb}[5]{
  {\begin{pmatrix*}[r]
    {#4} & {#3}  &  &  &  \\
    {#1} & {#2} & {#3} &  & \\
     & \ddots & \ddots & \ddots &  \\
     &  &  {#1} & {#2} &  {#3}\\
     &  &  & {#1} &  {#5} 
  \end{pmatrix*}}
}

\newcommand{\Uvec}{\mathcal{U}}

\newcommand{\Svec}{\mathcal{S}}
\newcommand{\Ld}{\mathcal{L}_\Delta}
\newcommand{\Nd}{\mathcal{N}_\Delta}

\newcolumntype{H}{>{\setbox0=\hbox\bgroup}c<{\egroup}@{}}